\documentclass[reqno]{conm-p-l}

\copyrightinfo{ }{ }

\theoremstyle{definition}

\theoremstyle{remark}

\numberwithin{equation}{section}

\usepackage[utf8]{inputenc}
\usepackage{wasysym}
\usepackage[colorinlistoftodos]{todonotes}
\usepackage{amsthm}
\usepackage{amssymb}
\usepackage{enumitem} 
\usepackage{mathrsfs} 
\usepackage{graphicx} 

\usepackage[unicode=true,pdfusetitle,
 bookmarks=true,bookmarksnumbered=false,bookmarksopen=false,
 breaklinks=false,pdfborder={0 0 1},backref=false,colorlinks=false]
 {hyperref} 

\theoremstyle{plain}
\newtheorem{thm}{Theorem}
\theoremstyle{definition}
\newtheorem{defn}[thm]{Definition}

\theoremstyle{remark}
\newtheorem{rem}[thm]{Remark}
\theoremstyle{remark}
\newtheorem*{rem*}{Remark}
\theoremstyle{plain}

\global\long\def\with{\,\middle|\,}
\def\Lsp{{\boldsymbol L}}
\global\long\def\d{\,{\rm d}}
\global\long\def\diff{{\rm d}}
\global\long\def\vertiii#1{\left\vert \kern-0.25ex  \left\vert \kern-0.25ex  \left\vert #1\right\vert \kern-0.25ex  \right\vert \kern-0.25ex  \right\vert }
\global\long\def\supp{\operatorname{supp}}

\newcommand{\N}{\setN}
\newcommand{\Z}{\setZ}

\newcommand{\R}{\setR}

\let\emptyset\varnothing

\def\Hiwsp{{\Hilb^1_{w}}}
\newcommand{\BMOsp}{{\boldsymbol {BMO}}}     
\newcommand{\Bsp}{{\boldsymbol B}}     
\newcommand{\BspN}{(\Bsp, \, \|\ebbes\|_\Bsp)}     
\newcommand{\ebbes}{\mbox{$\,\cdot\,$}}     
\newcommand{\Bspq}{{\Bsp^s_{\negthinspace p,q}}}     
\newcommand{\BspqHom}{{\dot{\Bsp}^s_{p,q}}} 
\newcommand{\BspqRd}{{\Bspqsp(\Rst^d)}}     
\newcommand{\Bspqsp}{{\Bsp^{s}_{\negthinspace p,q}}}     
\newcommand{\Rst}{{\mathbb R}}     
\newcommand\BspqRdN{{ (\BspqRd, \| \ebbes \|_{\Bspq}) }}     
\newcommand\Compl{\Cst} 
\newcommand{\Cst}{{\mathbb C}}     
\newcommand{\CooY}{{\Coosp(\Ysp)}}     
\newcommand{\Coosp}{{\boldsymbol{\mathcal Co}}}     
\newcommand{\Ysp}{{\boldsymbol Y}}     
\newcommand{\calD}{\mathcal{D}}     
\newcommand{\calQ}{\mathcal{Q}}     
\newcommand\FLi{{\mathcal F}{\negthinspace \Lisp}}     
\newcommand{\Lisp}{{\Lsp^1}}     
\newcommand{\FLiRd}{{ \FLi(\Rdst) }}     
\newcommand{\Rdst}{{{\Rst^d}}}     
\newcommand\FLiRdN{\big( \FLiRd, \, \|\ebbes\|_{\FLisp} \big)}     
\newcommand{\FLisp}{{ {\mathcal F}\Lisp}}     
\newcommand{\Fsp}{{\boldsymbol F}} 
\newcommand{\Fspq}{{\Fsp^s_{\negthinspace p,q}}}     
\newcommand{\FspqHom}{{\dot{\Fsp}^s_{p,q}}} 
\newcommand{\FspqRd}{{\Fspqsp(\Rst^d)}}     
\newcommand{\Fspqsp}{{\Fsp^{s}_{\negthinspace p,q}}}     
\newcommand\FspqRdN{{ (\FspqRd, \| \ebbes \|_{\Fspq}) }}     
\newcommand{\GL}{\operatorname{GL}}
\newcommand{\Hd}{{{\mathbb H}^d}}     
\newcommand\HiRd{{\Hsp^1 \negthinspace (\Rdst)}}     
\newcommand{\Hsp}{{\boldsymbol H}}     
\newcommand\Hilb{\mathcal H}     
\newcommand{\Liwsp}{{\Lsp^1_{\negthinspace w}}}     
\newcommand{\LtRd}{{\Ltsp(\Rst^d)}}     
\newcommand{\Ltsp}{{\Lsp^2}}     
\newcommand{\LtRdN}{\big( \LtRd, \, \|\ebbes\|_2 \big)}     
\newcommand{\LtcG}{{\Ltsp(\cG)}}     
\newcommand{\cG}{\mathscr{G}}     
\newcommand{\LtcGN}{\big( \LtcG, \, \|\ebbes\|_2 \big)}     
\newcommand\Msapq{{\Msp^{s,\alpha}_{\negthinspace p,q}}}     
\newcommand{\Msp}{{\boldsymbol M}}     
\newcommand{\MspqRd}{{\Mspqsp(\Rst^d)}}     
\newcommand\Mspqsp{{\Msp^s_{\negthinspace p,q}}}     
\newcommand\MspqRdN{\big( \MspqRd, \, \|\ebbes\|_{\Mspqsp} \big)} 
\newcommand{\ScRd}{{\Scsp(\Rst^d)}}     
\newcommand{\Scsp}{{\boldsymbol{\mathcal S}}}     
\newcommand\ScRdp{ { \Scsp}'(\Rdst)}     
\newcommand{\SchRd}{\ScRd}     
\newcommand{\SchpRd}{{\ScRdp}}     
\newcommand{\Wsp}{{\boldsymbol W}}     
\newcommand{\Xsp}{{\boldsymbol X}}
\newcommand{\YspN}{(\Ysp, \, \|\ebbes\|_\Ysp)}
\newcommand\Yspd{{\Ysp_{\nnnth d}}}     
\newcommand\nnnth{{ \negthinspace \negthinspace \negthinspace }}     
\newcommand{\bspq}{{{\boldsymbol{b}^s_{p,q}}}}     
\newcommand{\calA}{\mathcal{A}} 
\newcommand{\calB}{\mathcal{B}} 
\newcommand{\calH}{\mathcal{H}} 
\newcommand{\calO}{\mathcal{O}}     
\newcommand{\calP}{\mathcal{P}} 
\newcommand{\calR}{\mathcal{R}} 
\newcommand{\calS}{\mathcal{S}}     
\newcommand{\calU}{\mathcal{U}} 
\newcommand\coeff{{\mathcal C}}     
\newcommand{\decompsp}[3]{\calD ( {#1}, \Lsp^{#2}, {#3})}
\newcommand{\fsn}[2]{{\norm{#1}_{#2}}}     
\newcommand{\norm}[1]{\lVert#1\rVert}     
\newcommand{\fspq}{{\boldsymbol{f} ^s_{p,q}}}     
\newcommand{\lsp}{{\boldsymbol\ell}}     
\newcommand{\ltsp}{{\lsp^2}}     
\newcommand\ltspN{{\big( \ltsp, \, \|\ebbes\|_2 \big)}}     

\newcommand\recon{{\mathcal R}}     
\newcommand\suth{{ \, | \, } }     
\newcommand\wst{ w^{*} }     
\newcommand{\setR}{\mathbb{R}}
\newcommand{\setZ}{\mathbb{Z}}
\def\HiwA{{(\Hilb^1_{w})^{\angle}}}
\def\NSP#1{{(#1,\| \ebbes \|_{#1})}}
\newcommand{\setN}{\mathbb{N}}
\makeatletter
\def\@setcopyright{\relax}
\makeatother

\begin{document}

\title[From Frazier-Jawerth characterizations to Wavelets and Decomposition Spaces]{From Frazier-Jawerth characterizations of Besov spaces \\ to Wavelets and
Decomposition spaces}


\author{H.\@ G.\@ Feichtinger}
\address{Faculty of Mathematics, University of Vienna}  
\email{hans.feichtinger@univie.ac.at}

\author{F.\@ Voigtlaender}
\address{Lehrstuhl A für Mathematik, RWTH Aachen University}
\email{felix.voigtlaender@rwth-aachen.de}


\date{March 14, 2016} 

\begin{abstract}
This article describes how the ideas promoted
by the fundamental  
papers published by M.~Frazier and B.~Jawerth
in the eighties 
have influenced subsequent developments related to the
theory of {\it atomic decompositions} and {\it Banach frames}
for function spaces such as the {\it modulation spaces} and
{\it Besov-Triebel-Lizorkin spaces}.

Both of these classes of spaces arise as special cases of
two different, general constructions of function spaces:
{\it coorbit spaces} and {\it decomposition spaces}.
Coorbit spaces are defined by imposing certain decay conditions
on the so-called {\it voice transform}
of the function/distribution under consideration.
As a concrete example, one might think of the wavelet transform,
leading to the theory of Besov-Triebel-Lizorkin spaces.

{\it Decomposition spaces}, on the other hand, are defined using certain
decompositions in the Fourier domain. For Besov-Triebel-Lizorkin spaces, one uses a dyadic decomposition,
while a uniform decomposition yields modulation spaces.

Only recently, the second author has established a fruitful
connection between modern variants of wavelet theory with respect to
general dilation groups (which can be treated in the context of coorbit
theory) and a particular family of decomposition spaces. In this way,
optimal inclusion results and invariance properties for a variety of
smoothness spaces can be established. We will present an outline of these
connections and comment on the basic results arising in this context.
\end{abstract}

\maketitle

\section{The Frazier-Jawerth Atomic Characterizations}
\label{sec:FrazierJawertAtomicCharacterizations}

\subsection{Introductory statement}
This article is part of a series of notes (see e.g. \cite{fe14,fe15,fe09}),
which describe the role of different function spaces,
their various characterizations and their possible
applications from a ``postmodern viewpoint'',
emphasizing the importance of the concept of Banach frames (first defined in \cite{gr91}). For this development, a series of papers by Michael Frazier
and Bjoern Jawerth (\cite{frja85,frja88,frja90}) have been of great relevance. Hence
we will try to demonstrate how the ideas of these papers lived on and expanded
in the subsequent decades.
\pagebreak

\subsection{Notation}
\label{sub:Notation}
We will employ the following notation:
For a group $G$ and any function $f : G \to S$, for some set $S$, we define
\begin{align*}
    L_x f : G \to S, y \mapsto f(x^{-1}y), 
    \quad
    R_x f : G \to S, y \mapsto f(yx), \quad  x,y \in G.
\end{align*}
In the special case of the (abelian) group $\R^d$, we also write $T_x := L_x$. Furthermore, for $f : \R^d \to \Compl$ and $\xi \in \R^d$, we define the
\textbf{modulation} of $f$ by $\xi$ as
\[
    M_\xi f : \R^d \to \Compl, x \mapsto e^{2\pi i x \cdot \xi} \cdot f(x).
\]
Finally, for $h \in \GL (\R^d)$, we define the ($\Lsp^2$ normalized) \textbf{dilation} of $f$ by $h$ as
\[
    D_h f : \R^d \to \Compl, x \mapsto |\det h|^{-1/2} \cdot f(h^{-1} x).
\]
For the special case $h = a \cdot {\mathrm{id}}$ with $a \in \R \setminus \{0\}$, we also write $D_a := D_{a \cdot {\mathrm id}}$.

For the \textbf{Fourier transform}, we use the normalization
\[
    \mathcal{F}f ( \xi) := \widehat{f}(\xi) := \int_{\R^d} f(x) e^{-2\pi i x\cdot \xi} \diff x, \quad f \in \Lsp^1(\R^d).
\]
It is well-known that the Fourier transform extends to a unitary automorphism of $\Lsp^2(\R^d)$, where the inverse is the unique extension of
the operator $\mathcal{F}^{-1}$, given  by
\[
    \mathcal{F}^{-1} h ( x) = \int_{\R^d} h(\xi) e^{2\pi i x \cdot \xi} \diff \xi,
    \quad h \in \Lsp^1(\R^d).
\]
$\lambda (M)$ denotes  the $d$-dimensional Lebesgue measure of a (measurable) set $M \subset \R^d$.

Finally, for $x \in \R$, we write $x_+ := \max\{x,0\}$.

\subsection{The Essence of the work of Frazier-Jawerth}

To the best of our knowledge, the influential papers \cite{frja85,frja88,frja90}
have been the first to fully characterize two families of Banach spaces of tempered
distributions, namely the Besov spaces $\BspqRdN$ and the Triebel-Lizorkin
spaces  $\FspqRdN$, by corresponding growth and summability conditions on a
suitably defined sequence of coefficients. These coefficients depend linearly on
the function/distribution under consideration. Similar atomic representation
theorems had been realized only a few years earlier in the context of harmonic
function spaces (see e.g.\ \cite{coro80,rita83})

In a more modern terminology (going back to the work of K.~Gr\"ochenig \cite{gr91}),
one could say that Frazier and Jawerth established specific, but rather concrete {\it Banach frame} expansions
for these two families of function spaces, starting from the characterization of these spaces via
{\it dyadic partitions of unity} on  the Fourier transform side.
This characterization in turn is based on the description of Besov-Triebel-Lizorkin spaces via dyadic decompositions
(see \cite{pe76,tr77,tr83}).

In this general theory of Banach frames, one has---roughly speaking---the following situation:

(i) Starting from a fixed family of (possibly non-orthogonal, but sufficiently
      rich) {\it atoms} $(g_i)_{i \in I}$, one obtains \emph{atomic representations} of the form $f =  \sum_{i \in I} c_i g_i$, typically with convergence in the norm of $\BspN$ (or at least $\wst$-convergence for the case of dual spaces).

(ii) The coefficients are obtained in a linear way, i.e.,
     there are suitable bounded linear {\it coefficient mappings}
     $f \mapsto c_i = c_i (f)$ satisfying
     $f = \sum_{i \in I} c_i g_i$. These mappings are often realized as scalar
     products with respect to a suitable family of atoms $(h_i)_{i \in I}$, forming a so-called {\it dual frame}.
     Thus, $c_i = \langle f, h_i \rangle$, where the ``scalar product'' can be viewed as an extension of the Hilbert-space scalar product.

(iii) Although the representation is by far not unique, there is a high
degree of compatibility between the membership of $f$ in one of the spaces
under consideration and membership of the coefficients $(c_i)_{i \in I}$ in
a suitable Banach space of sequences, in fact, in a suitable {\it solid BK space}.
These spaces have the property that sequences which are smaller in terms of absolute values---in a coordinate-wise sense---also have a smaller norm.

Precisely, it is claimed that for a function space $\Xsp$ (from a certain class),
one can find a BK space $\Ysp \leq \Compl^I$, a bounded {\it coefficient map}
$\coeff : \Xsp \to \Ysp, f \mapsto (c_i(f))_{i \in I}$,
and a {\it synthesis mapping} $\recon : \Ysp \to \Xsp, {\bf c} \mapsto \sum_{i \in I} c_i g_i$ which is a left inverse for $\coeff$.\\[-0.1cm]

More concretely, the \textbf{$\varphi$-transform} as introduced by Frazier and Jawerth in \cite{frja85,frja88,frja90} is a \emph{fixed, linear} transformation
    $S_\varphi : f \mapsto \left((S_\varphi f)_Q \right)_{Q \in \mathcal{Q}}$
with the benefit that a variety of function spaces can be characterized in terms of the \emph{size} of $S_\varphi f$, i.e.\@ by
$\fsn {S_\varphi f} \Ysp$ for suitable \emph{solid} sequence spaces $\Ysp$. Solidity of $\Ysp$ formally encodes the requirement that only the \emph{size}
of the coefficients should be important, so that there is no cancellation between the different coefficients.

The class of spaces for which a characterization in terms of the $\varphi$-transform is possible includes the classes of \emph{Besov spaces} $\Bspq$ and the
\emph{Triebel-Lizorkin spaces} $\Fspq$. Note that both of these classes of spaces come in two variants: homogeneous and inhomogeneous spaces.
The theory developed by Frazier and Jawerth applies to both of these subclasses. For concreteness, we will concentrate in the sequel
on the \emph{in}homogeneous spaces,
which have the advantage of being modulation invariant.

To describe the $\varphi$-transform and the resulting decomposition results more precisely, we begin with the index set $\mathcal{Q}$, which is the set
of all \textbf{dyadic cubes} with side-length $\leq 1$. Here, by definition (cf.\@ \cite{frja90}), a dyadic cube is a set of the form
\[
    Q_{\nu, k} = \left\{ x \in \R^d \with \forall i \in \{1, \dots, d\}: 2^{-\nu} k_i \leq x_i < 2^{-\nu} (k_i + 1) \right\}
\]
for arbitrary $\nu \in \Z$ and $k \in \Z^d$.
Given such a cube $Q = Q_{\nu, k}$, we define the \textbf{lower left corner} of $Q$ as $x_Q := 2^{-\nu}k$
and the \textbf{side length} of $Q$ as $\ell(Q) := 2^{-\nu}$.
To distinguish cubes with different side lengths, we also define for $\nu \in \Z$ the set
\[
    \mathcal{Q}_\nu := \left\{ Q \in \mathcal{Q} \with \ell(Q) = 2^{-\nu} \right\}.
\]
Finally, for an arbitrary function $\varphi : \R^d \to \Compl$ and $Q = Q_{\nu, k}$, we let
\[
    \varphi_Q (x) := 2^{\nu d /2} \cdot \varphi(2^{\nu} x - k) = 2^{-\nu d/2} \cdot \varphi_{\nu} (x - x_Q),
\]
where $\varphi_\nu (x) := 2^{\nu d} \cdot \varphi (2^{\nu} x)$. Note that if $\varphi$ is ``concentrated'' in $[0,1]^d$, then $\varphi_Q$ is ``concentrated'' on $Q$.
Furthermore, we have $\fsn {\varphi_Q} {\Lsp^2} = \fsn \varphi {\Lsp^2}$.

Now, the $\varphi$-transform is defined using two \textbf{analyzing windows} $\varphi, \varphi^0 \in \SchRd$.
 The basic assumption regarding these windows
is that there are ``dual windows'' $\psi, \psi^0 \in \SchRd$
satisfying the following (cf.\@ \cite[eqs. (12.1)--(12.2) and (2.1)--(2.3)]{frja90})
\begin{align}
    \supp \widehat{\phi^0}, \widehat{\psi^0} &\subset
       \left\{ \xi \in \R^d \with |\xi| \leq 2 \right\}, \\
    \supp \widehat{\varphi}, \widehat{\psi} &\subset \left\{ \xi \in \R^d \with 1/2 \leq  |\xi| \leq 2 \right\}, \\
    \overline{\widehat{\varphi^0}(\xi)} \cdot \widehat{\psi^0}(\xi) + \sum_{\nu = 1}^{\infty} \overline{\widehat{\varphi}(2^{-\nu}\xi)} \cdot \widehat{\psi}(2^{-\nu}\xi) &=1  \qquad \forall \xi \in \R^d.
\end{align}
Under this assumption, the \textbf{(inhomogeneous) $\varphi$-transform} $S_\varphi$ is the analysis operator
of the frame, so that the $\varphi$-transform $S_\varphi f = \left( (S_\varphi f)_Q \right)_{Q \in \mathcal{Q}}$ of a
tempered distribution $f \in \SchpRd$ is defined by
\[
    (S_\varphi f)_Q := \begin{cases}
                            \langle f, \varphi_Q\rangle ,& \text{if } \ell(Q) < 1,\\
                            \langle f, \varphi^0_Q\rangle    ,& \text{if } \ell(Q) = 1.
                       \end{cases}
\]
Finally, the (formal) inverse of $S_\varphi$---the corresponding synthesis
operator---is
\[
    T_\psi [(s_Q)_{Q \in \mathcal{Q}}] := \sum_{Q \in \mathcal{Q}_0} s_Q \psi^0_Q  + \sum_{\nu = 1}^{\infty} \sum_{Q \in \mathcal{Q}_{\nu}} s_Q \psi_Q.
\]
Now, as shown in \cite[equation (12.4)]{frja90}, we have $f = T_\psi S_\varphi f$ for all $f \in \SchpRd$, with convergence
of the series defining $T_\psi (S_\varphi f)$ in the topology of $\SchpRd$.

The main result of Frazier and Jawerth concerning inhomogeneous Triebel-Lizorkin spaces using the $\varphi$-transform reads as follows (cf.\@ \cite[Theorem 12.2]{frja90}):
\begin{thm}
    \label{thm:FrazierJawerthInhomogeneousFspqCharacterization}
    For $s \in \R$, $0<p<\infty$ and $0<q\leq \infty$, the operators
    \[
        S_\varphi : \Fspq \to \fspq \qquad \text{ and } \qquad T_\psi : \fspq \to \Fspq
    \]
    are bounded, with $T_\psi \circ S_\varphi = \mathrm{id}$ on $\Fspq$. Hence,
    $\fsn f \Fspq \asymp \fsn {S_\varphi f} \fspq$ and $\Fspq$ is a retract of $\fspq$, i.e., it can be identified with a complemented subspace
    of $\fspq$.
\end{thm}
In the above theorem, the solid BK-space $\fspq$ is
defined as the set of all sequences
$c = (c_Q)_{Q \in \mathcal{Q}} \in \Compl^{\mathcal{Q}}$ for which
the (quasi)-norm
\begin{align*}
\fsn c \fspq & := \left\Vert \left\Vert \left(\left[\ell\left(Q\right)\right]^{-s}\cdot\left|c_{Q}\right|\cdot\widetilde{\chi_{Q}}\right)_{Q\in\mathcal{Q}}\right\Vert _{\lsp^{q}\left(\mathcal{Q}\right)}\right\Vert _{\Lsp^{p}\left(\d x\right)} 
\end{align*}
is finite. 
Here, we used the symbol
$
    \widetilde{\chi_E} := (\lambda(E))^{-1/2} \cdot \chi_E
$
for the $\Lsp^2$-normalized version of $\chi_E$, for any measurable set $E \subset \R^d$ of finite, positive measure.

In addition to the above result for Triebel-Lizorkin spaces, \cite{frja85} provides the following analogue characterization of Besov spaces.
\begin{thm}
    \label{thm:FrazierJawerthInhomogeneousBesovCharacterization}
    For $0 < p,q \leq \infty$ and $s \in \R$, the operators
    \[
        S_\varphi : \Bspq \to \bspq \qquad \text{ and } \qquad T_\psi : \bspq \to \Bspq
    \]
    are bounded, with $T_\psi \circ S_\varphi = \mathrm{id}$ on $\Bspq$.
    Hence, $\fsn f \Bspq \asymp \fsn {S_\varphi f} \bspq$ and $\Bspq$ is a
    {\it retract} of $\bspq$, and can be identified with a complemented subspace
    of $\bspq$.
\end{thm}
The solid BK-space $\bspq$ consists of all $c = (c_Q)_{Q \in \mathcal{Q}} \in \Compl^{\mathcal{Q}}$ for which the following
quasi-norm is finite:
\begin{align*}
    \fsn c \bspq &:= \left\Vert \left(\left\Vert \left(\left[\ell\left(Q\right)\right]^{d\left(\frac{1}{p}-\frac{1}{2}\right) - s} c_{Q}\right)_{Q\in\mathcal{Q}_{\nu}}\right\Vert _{\lsp^{p}(\mathcal{Q}_\nu)}\right)_{\nu\in\mathbb{N}_{0}}\right\Vert _{\lsp^{q}\left(\mathbb{N}_{0}\right)}. 
\end{align*}
Similar results have later been described in the books of H.~Triebel,
who showed that multivariate wavelet orthonormal bases are not just
bases for the Hilbert space $\LtRd$, but also (in a modern terminology)
{\it Riesz projection bases} for the corresponding solid BK-spaces.
This is another way of saying that the unitary isomorphism
between $\Lsp^2 (\R^d)$ and $\lsp^2(I)$---which is induced by the orthonormal basis $(\varphi_i)_{i \in I}$---extends
to an isomorphism between Banach spaces and their discrete counterparts,
for a whole family of spaces and with uniform bounds.
But the existence of orthonormal wavelet bases was established (in \cite{leme86,me88,da88-1}) only after
the appearance of the paper \cite{frja85} of Frazier and Jawerth.

In summary, these ``decomposition'' results of Frazier and Jawerth imply
a variety of statements:  \\
(1)    A \emph{consistency} statement: Any two (admissible)
        window-families $\varphi_1, \varphi^0_1$ and $\varphi_2, \varphi^0_2$
        yield the \emph{same} space by imposing certain decay conditions
        on the associated $\varphi$-transform, i.e., we have for $ i = 1,2$:
        \[
            \left\{f \in \SchpRd \with S_{\varphi_i} f\in \fspq \right\} = \Fspq .
        \]
(2)   A \emph{decomposition} (or \emph{representation}) statement:
        \[
            f = \sum_{Q \in \mathcal{Q}} s_Q (f) \cdot \psi_Q,
            \quad f \in \Fspq,
        \]
        with $(s_Q (f))_{Q \in \mathcal{Q}} = S_\varphi f \in \fspq$
        and atoms $\psi_Q$ of the special form
        \[
            \qquad \qquad \psi_Q (x) = 2^{-\nu d/2} \cdot \psi_{\nu} (x - x_Q) = [\pi(x_Q, 2^{-\nu}) \psi](x) \quad \text{ for } \quad \ell(Q) = 2^{-\nu},
        \]
        where $\pi(x, a) := T_x D_a$ is the
        quasi-regular representation on $\Lsp^2 (\R^d)$ of the $ax+b$ group.
        This 
        was already observed in \cite[Remark 3.2]{frja90}.

(3)    The two parts of the representation 
       are all valid for a certain \emph{range} of spaces, not just for a
       \emph{single} Hilbert space or one \emph{individual} Banach space.\label{FV:SpecialTarget}

        It is also worth noting that, if one restricts the construction of Frazier and Jawerth to the setting of Hilbert spaces,
        one finds that they are implicitly establishing that the set of
        analyzing atoms forms a frame; but the elements used to perform
        the reconstruction are not the elements of the (canonical) dual frame.


(4)    The concrete BK-spaces $\fspq$ arising in this context are not just
        abstract BK-spaces (Banach spaces of numerical sequences, whose
        coordinates depend continuously on the sequence). In addition, they are
        {\it Banach lattices} (following the terminology of Luxemburg-Zaanen)
        resp.\ so-called {\it solid} BK-spaces: With each sequence
        ${\bf d} = (d_Q)_{Q \in \calQ} \in \fspq$, also any sequence $ {\bf b} = (b_Q)_{Q \in \calQ}$ with
        $ |b_Q| \leq |d_Q|$ for all $Q \in \calQ$ belongs to $\fspq$, with
        $ \|\bf b\|_{\fspq} \leq \|\bf d\|_{\fspq}$.\\[-0.1cm]

The same statements remain true with $\Bspq$ and $\bspq$ instead of $\Fspq$ and $\fspq$.

The above observations are similar to the main properties of \emph{coorbit spaces}, as introduced by Feichtinger and Gröchenig in \cite{fegr88,fegr89,fegr89-1}, see also Section \ref{sec:CoorbitTheory}.


\begin{rem*}
While the $\varphi$-transform and the resulting atomic decompositions for
Besov (resp.\ Triebel-Lizorkin) spaces grew out of the characterization
of these spaces using Littlewood-Paley theory---essentially established
by the earlier work of J.~Peetre and  H.~Triebel---it turned out that
they were also fore-runners for the characterization of these function
spaces via {\it wavelet frames}.
\end{rem*}

\subsection{Modulation Spaces}  
\label{sec:ModulationSpace}

One of the earliest variations of the atomic decomposition method
proposed by Frazier and Jawerth arose in connection
with so-called {\it modulation spaces}.
These spaces have been introduced with the idea of describing smoothness for
functions over locally compact Abelian groups (LCA groups) 
via {\it uniform} (instead of dyadic) decompositions on the Fourier transform
side. The first report \cite{fe83-4} issued in 1983 (see also \cite{fe03-1} for an expanded form)
came too early to be appreciated at the time, certainly also because
it was formulated in the most general setting.
A large number of references concerning these spaces can be found in \cite{fe06}.
Nowadays, modulation spaces are viewed as natural domains for time-frequency analysis and certain families of pseudo-differential operators.

Although designed from the very beginning at the most general level,
a specific subfamily of modulation spaces, namely the modulation spaces
$\MspqRdN$, with $1 \leq p,q \leq \infty, s \in \Rst$,
are most similar to the corresponding
family of Besov spaces $\BspqRdN$, with the same parameters. In fact (only) for the
case $p=q=2$ these spaces coincide, while otherwise they are different (cf.\@ the
PhD thesis of P.~Gr\"obner, \cite{gr92-2}).
By now, there is an extensive literature concerning the inclusion relation
between different types of modulation spaces resp.\ more general decomposition
spaces, see e.g. \cite{ok04,to04-1,suto07,hawa14,VoigtlaenderPhDThesis}.

A Frazier-Jawerth type atomic characterization of modulation spaces
has first been given in \cite{fe89-1} (received by the editors in Oct.\@ 1986).
These results are quite similar to those of Frazier-Jawerth,
although no direct reference is made to their papers, obviously because
their results have not been used. Instead, the argument is based on
 a combination of Shannon's theorem
and techniques from the theory of Wiener amalgam spaces developed in the
early 80s (\cite{fe83,fe92}). These techniques provide refined variants of Poisson's formula,
which are used to obtain what is nowadays called a {\it Gabor characterization}
of modulation spaces.

To emphasize the analogy between \cite{frja85,frja88} and \cite{fe89-1},
we note the following: In each case, suitable
partitions of unity are used, which are
bounded families within the Fourier algebra $\FLiRdN$.
Further, the building blocks are carefully chosen by the authors in a specific way, to
achieve the goal of atomic representation.

In contrast to the last point, the typical question treated e.g.\@ in the context
of Gabor analysis---or more generally coorbit theory---is the following:
{\it Given}  some (structured) family of atoms (so without chance of the
user to match the atoms with the function spaces), can one still verify
frame properties and compute a dual frame, or at least guarantee the existence of such a dual frame?

Compared to the situation in wavelet theory, one can say:
The possibility of generating {\it orthonormal wavelet bases} came certainly
as a surprise, especially to Yves Meyer, who tried to prove the converse,
but ended up with his first construction of such a basis.  His negative
expectations were probably based on his knowledge of the Balian-Low Theorem
which excludes the existence of Riesz bases (and thus also of orthonormal bases) of Gaborian type with ``good atoms''
which are well concentrated in the time-frequency sense (\cite{ba81-1,behewa95}).
The rapid development of wavelet theory, starting with \cite{leme86}, is in large parts due to the possibility of having {\it orthonormal wavelet bases} which may be difficult to construct, but which are easy to use. The construction of compactly
supported orthonormal wavelet bases with a given amount of smoothness
by Ingrid Daubechies (\cite{da88-1})
was one of the main reasons why the Frazier-Jawerth decomposition method
was not pursued too much for a while.

A crucial property of both orthonormal wavelets
as well as the Schwartz building blocks for the Frazier-Jawerth decompositions
was the fact that from the very beginning, the canonical isomorphism between the Hilbert space $\LtRdN$ and $\ltspN$ extends
automatically to an isomorphism between the classical function spaces
(those from the family of Besov-Triebel-Lizorkin type, including $\HiRd$ and
$\BMOsp$) and (subspaces of) the corresponding solid BK-spaces (typically weighted mixed-norm sequence spaces).

It is interesting to note that the sequence spaces used for the Frazier-Jawerth
characterizations and the orthonormal wavelet systems are more or less the same,
once one identifies the collection of cubes with different centers with the unique
affine transformation (from the $ax+b$-group) needed to obtain it from a symmetric standard cube.
The common structure of all the different types of atomic characterizations
going back to Frazier-Jawerth is the
identification of a family of Banach spaces (of functions or distributions) via some coefficient mapping with suitable closed and complemented subspace of the corresponding family of solid BK-spaces. 

\section{Coorbit Theory}
\label{sec:CoorbitTheory}

\subsection{Unifying Approach to the CWT and the STFT}
\label{sec:UnifCoorbitTheory}

Once the theory of wavelets started to gain momentum,
with the classical function spaces being characterized
via the Frazier-Jawerth decompositions or via suitable
(orthonormal) wavelet systems 
and once 
it was understood that there are
similar characterizations for modulation spaces
via \emph{Gabor expansions}, it became a natural question to ask
whether these two signal representation methods have anything further in common.

Group representation theory finally provided such a link:
in both (as well as other) cases there is a  {\it (square) integrable}
group representation of some locally compact group, acting in an irreducible way
(more or less) on a given Hilbert space. In the language of group representation
theory, the relevant group for the theory of wavelets is the {\it affine group} (also called ``$ax+b$''-group), acting by dilations and translations on $\Lsp^2 (\R^d)$. On the other hand, for the theory of modulation spaces or Gabor expansions, the group acting in the background is the {\it (reduced) Heisenberg group}
$\Hd$, acting on $\Lsp^2 (\R^d)$ by translations and modulations.

\subsection{The Essence of Coorbit Theory}
\label{sec:EssCoorbitTheory}

The core of {\it coorbit theory}, as outlined in \cite{fegr89,fegr89-1},
is to generate---from a given unitary group representation $\pi$ of some
locally compact group $\cG$ on a Hilbert space $\Hilb$---a whole
family of spaces $\CooY$, the so-called {\it coorbit spaces} (for $(\cG,\pi)$).
For large subfamilies of these spaces, coorbit theory provides atomic decompositions which allow to characterize elements of the space $\CooY$ by the fact that they can be described using coefficients in the solid BK-spaces  $\Ysp_d$ which is naturally associated with $\CooY$.

The first step for describing $\CooY$ is to define the {\it voice transform} $V_g$,
by
$$(V_g f)(x) := \langle f,\pi(x) g\rangle, \qquad x \in \cG, \quad f,g \in \Hilb.$$
If $g$ is suitably chosen (i.e.\@ satisfying $\| V_g g \|_{\Lsp^2} = \|g\|_{\Hilb}$),
the mapping $f \mapsto V_g f$ maps $\Hilb$ isometrically\cite{grmopa85} into $\LtcGN$, so that $\Hilb_0 := V_g(\Hilb)$
is a closed, left-invariant subspace of $\LtcG$.
Note that this crucially uses that the representation $\pi$ is (square) integrable and irreducible.

The name ``voice transform'' goes back to the paper \cite{gogrmo84}, where this name is used
synonymously with the term ``cycle-octave transform'' for what later
became the {\it continuous wavelet transform} (CWT). Since
in that paper both the Gabor and the wavelet case have been addressed,
this terminology was taken into the general coorbit theory. It was also
used later on in the context of the Blaschke group by Margit Pap and
Ferenc Schipp (see e.g. \cite{pasc06,pa10-1,pa12}).

Given the voice transform $V_g$ and a translation invariant, solid BF-space%
\footnote{This means that $\YspN$ is a space of measurable functions on $\cG$ such that if $f \in \Ysp$ and
    if $g : \cG \to \Compl$ is measurable with $|g|\leq |f|$ almost everywhere, then $g \in \Ysp$ with $\|g \|_{\Ysp} \leq \|f \|_{\Ysp}$.}
    $\YspN$ on $\cG$, we \emph{would like} to define the associated coorbit space $\CooY$ as
$\CooY := \{ f \in \Hilb \suth V_g f \in \Ysp \} $, with the norm $\fsn f \CooY := \fsn {V_g f} \Ysp$.
But this will in general \emph{not} yield a complete space.

Recalling that the classical function spaces (like Besov spaces) are defined as subsets of the space $\SchpRd$ of tempered distributions
and not of $\Lsp^2 (\R^d)$, we thus have to find a suitable replacement for the class of Schwartz functions in the present generality.
In this setting, however, there will (in general) be no ``universal'' class of test functions like the Schwartz space. Instead, given the solid BF-space $\YspN$,
one has to choose a suitable class of test functions \emph{specific} to $\YspN$.

For this, one first needs to choose a so-called \emph{control weight}
$w : \cG \to (0,\infty)$ for the space $\YspN$.
This weight depends in a certain way on the operator norms of left- and right translation on $\YspN$.
 For the sake of brevity we omit these details.
 Note, however, that $w$ is assumed to be submultiplicative
 (i.e., $w(xy) \leq w(x) \cdot w(y)$).
Given such a control weight $w$,  we introduce the space
\begin{equation}
    \Hiwsp := \{ f \in \Hilb \suth V_g f \in \Liwsp\}, \qquad \text{ where} \qquad \Liwsp = \{ f \suth fw \in \Lisp(\cG) \},
    \label{eq:CoorbitTheoryTestFunctions}
\end{equation}
which will play the role of ``test functions'' in the present setting; the coorbit spaces will thus be subspaces of the (anti)dual space $\HiwA$ of $\Hiwsp$.

Of course, the ``analyzing window'' $g$---which is used to define the space $\Hiwsp$---needs to be chosen carefully.
As shown in \cite{fegr88,fegr89}, the space $\Hiwsp$ is independent of the precise choice of $g$, as long as
\begin{equation}
    0 \neq g \in \calA_w (\cG) := \left\{  g \in \Hilb \with V_g g \in \Liwsp \right\}.
    \label{eq:AnalyzingWindowClass}
\end{equation}
The space $\calA_w$ is called the space of \emph{analyzing windows}. A crucial \emph{assumption} for the applicability of coorbit theory
is that $\calA_w \supsetneq \{0\}$ is nontrivial.

Later on we will need the class of \emph{better vectors}
$\calB_w \subset \calA_w$, whose precise definition we omit.
For the initiated reader, we note that
$g \in \calB_w$ needs to satisfy $V_g g \in \Wsp^R (\Lsp^\infty, \Liwsp)$, where $\Wsp^R$ denotes a \emph{(right-sided) Wiener amalgam space}.

Before defining  {\it coorbit spaces} in full generality, we need to extend
the voice transform to the reservoir $\HiwA$:
Due to the $\pi$-invariance of $\Hiwsp$,
the (generalized) {\it voice transform}
$$ (V_g f)(x) = \langle f, \pi(x)g \rangle, \quad x \in \cG,$$
of $f \in \HiwA$ with respect to $g \in \Hiwsp$ is well-defined.
Here, the pairing is understood as $\langle f, g\rangle = f(g)$,
which extends the usual scalar product on $\Hilb$.

With all this, we define 
the corresponding {\it coorbit space}
$$ \CooY := \{ f \in \HiwA \suth V_g f \in \Ysp \},$$
with the usual natural norm. Again one can show that the space is \emph{independent} of the analyzing window $0 \neq g \in \calA_w$.
Another main result of coorbit theory concerns the atomic decomposition claim:
%
\begin{thm}
    \label{thm:CoorbitAtomicDecomposition}
For any $0 \neq g \in \calB_w$ there exists a neighborhood $U$ of the identity $e \in \cG$ such that the following is true:
For any {\it well-separated}\footnote{This means that $(x_i)_{i \in I}$ is the finite union of uniformly separated sets in $\cG$.}
family $(x_i)_{i \in I}$ in $\cG$ which is $U$-dense,
the family $(\pi(x_i) g)_{i \in I}$ defines a Banach frame for $\CooY$.

More precisely: There exists a solid BK-space $\Ysp_d = \Ysp_d ( (x_i)_{i \in I} ) \leq \Compl^I$
and a bounded linear mapping
$\coeff : f \mapsto (c_i (f))_{i \in I}$ from $\CooY$ to $\Ysp_d$,
such that 
$$ f = \recon (\coeff f) = \sum_{i \in I}  c_i (f) \pi(x_i)g, \quad  \forall \, f \in \CooY, $$
with (unconditional) convergence in  $\NSP \CooY$, if $\Hiwsp$
is dense in $\CooY$ (and in the $\wst$-sense otherwise). Furthermore, $\recon : \Ysp_d \to \CooY$ is bounded.
\end{thm}

In comparison with the original Frazier-Jawerth approach we note the following:

(a)
In the FJ-approach both the analyzing and the synthesizing vector
    are specified by the construction; although there is some freedom in the construction, the two elements have to be
    chosen in a very specific way.

(b)
In contrast, coorbit theory provides a much larger reservoir of atoms;
in fact, coorbit theory shows that one can reconstruct the \emph{complete} voice transform $V_g f$ (and hence the element $f \in \CooY$)
from its {\it samples} $(V_g f (x_i))_{i\in I}$, as long as $g \in \calB_w$ and as long as the samples are taken densely enough.

(c) Both theories provide a
retraction between a family of function spaces
and the corresponding solid BK-spaces. Since the
$ax+b$-group and the (reduced) Heisenberg group both have two unbounded
variables, it is natural to use mixed-norm conditions on the variables.

The corresponding discrete variants are discrete mixed-norm spaces,
as long as the sampling points $(x_i)_{i \in I}$ form a product set.
For a {\it rearrangement invariant} space $\Ysp$, the
sequence space $\Yspd$ is just the ``natural'' discrete variant; e.g.\ $(\Lsp^p)_d = \ell^p$.

(d) In both the Frazier-Jawerth-theory and the coorbit approach, one has the choice
to describe them in full generality or for a family of spaces.
The published versions of Frazier-Jawerth-theory describe the decompositions
in the setting of tempered distributions; consequently they are
able to characterize Besov spaces for all real parameters $s$.

If one relaxes the conditions on the building blocks
(e.g.\@ if one assumes only a certain number of vanishing moments
instead of requiring the building blocks to be compactly supported
away from the origin),
then one would have a valid statement only for a certain range of
spaces $\BspqRd$, with $|s| \leq s_0$.

In the coorbit setting, the family of spaces $\CooY$ for which a
certain window $g$ is suitable is determined by the weight $w$:
An analyzing window $g \in \calA_w$
is only guaranteed to characterize $\CooY$, i.e., to satisfy
$
    \CooY = \left\{ f \with V_g f \in \Ysp \right\},
$
if $w$ is a control weight for $\Ysp$. Likewise, Theorem \ref{thm:CoorbitAtomicDecomposition} only guarantees
that $g \in \calB_w$ generates a Banach frame for $\CooY$ if $w$ is a control weight for $\Ysp$.

\subsection{Shearlets and other constructions}
\label{sec:BeyCoorbitTheory}

The first aim of the theory of {\it coorbit spaces} as proposed
in \cite{fegr88,fegr89,fegr89-1,fegr92} was certainly to provide a unified treatment of the two most important situations (see Subsection \ref{sec:UnifCoorbitTheory})
where {\it painless non-orthogonal expansions} arose (\cite{dagrme86}).
But it also paved the way to consider other examples through the lens of the general theory.

Already early on, it was clear that certain {\it Moebius invariant}
Banach spaces of analytic functions on the unit disk
(see the work of Arazy-Fisher-Peetre \cite{arfipe85})
fit into the framework of coorbit theory.
These spaces are described by the behaviour of their members,
typically by imposing integrability properties, expressed by
weighted mixed norm spaces, with a radial and a circular component.
In a way, they can be compared with Fock spaces, which consist of
analytic functions over phase space, but can be identified with
the set of voice transforms (i.e.\ STFTs) with Gaussian window.

After shearlets were first introduced by Kutyniok, Labate, Lim, Guo and
Weiss\cite{kulaliwe05,gula07-1}, their group theoretic nature was realized in
\cite{dakumasastte08}, see also \cite{kula09}. Building upon that group-theoretic background of the shearlet transform,
the theory of {\it shearlet spaces} was investigated 
in \cite{dakustte09,dastte10,dastte12,dastte11,dahastte13,dahate12}. In \cite{dakustte09}, applicability of coorbit theory (i.e.\@ $\calB_w \neq \{0\}$) in dimension $d=2$ was established and the associated Banach frames (as provided by Theorem \ref{thm:CoorbitAtomicDecomposition}) were written down explicitly.
The generalization to higher dimensions (including the definition of a possible shearlet group in dimensions $d>2$) was obtained in \cite{dastte10}.
Another generalization of the shearlet group to higher dimensions was considered in \cite{dahate12}.
Moreover, the relation of shearlet coorbit spaces to more classical smoothness spaces (namely to (sums of) homogeneous Besov spaces)
was investigated in \cite{dastte11} for $d=2$ and in \cite{dahastte13} for higher dimensions. Related results will be discussed in Section \ref{sec:WaveletTheory}.

Yet another group, the so-called {\it Blaschke group}, is in the background
of a series of papers by M.~Pap and her coauthors
(\cite{pasc06,pa10-1,pa12,fepa13}).



\section{Decomposition Spaces}
\label{sec:DecompositionSpaces}

A natural starting point for the definition of decomposition spaces is the observation that
modulation spaces, as well as Besov spaces, can be described by imposing certain decay conditions
on the sequence of $\Lsp^p$ norms $\left( \|\mathcal{F}^{-1} (\varphi_i \widehat{g})\|_{\Lsp^p} \right)_{i\in I}$ for suitable families of
functions $(\varphi_i)_{i \in I}$. In the case of (inhomogeneous) Besov spaces, the $(\varphi_n)_{n \in \N_0}$ form a \emph{dyadic} partition of unity,
while for modulation spaces, a \emph{uniform} partition of unity $(\varphi_k)_{k \in \Z^d}$ is used.

By utilizing this observation, \emph{decomposition spaces}, as introduced by Feichtinger and Gröbner in \cite{fegr85,fe87}, provide
a common framework for Besov- and modulation spaces, as well as many other smoothness spaces occurring in harmonic analysis.

Indeed, we will see that decomposition spaces can be used
to describe the $\alpha$-modulation spaces---a family of spaces intermediate to Besov- and modulation spaces.
Furthermore, they provide an alternative view on a large class of wavelet type coorbit spaces.
We will see that the decomposition space viewpoint
makes many properties of these spaces transparent, while it is highly nontrivial (if not impossible)
to obtain these properties directly using the coorbit viewpoint.

\subsection{Basic properties of decomposition spaces}
\label{sub:DecompositionSpaces}

The original{\textemdash}very general{\textemdash}definition of decomposition spaces\cite{fegr85,fe87} starts with a covering $\mathcal{Q} = (Q_i)_{i \in I}$ of a given space $X$.
Most of the time, one may think of $X$ as (a subset of) the frequency space $\R^d$; but in principle, it could be some manifold or other domain.

Now, given a suitable function space $\Bsp$, the idea is to define the \textbf{decomposition space norm} of a function/distribution $f$
by measuring the \emph{local} behaviour of $f$
in the $\Bsp$-norm, i.e.\@ by localizing $f$ to the different sets $Q_i$
(using a suitable partition of unity $(\varphi_i)_{i \in I}$)
and by measuring the individual pieces using the $\Bsp$-norm.
The \emph{global} properties of this (generalized) sequence of $\Bsp$ norms are then
restricted using a suitable sequence space $\Ysp$.
Formally, we define
\begin{equation}
    \|f \|_{\mathcal{D}(\mathcal{Q}, \Bsp, \Ysp)}   := \left\| \left(  \|\varphi_i \cdot f\|_\Bsp \right)_{i \in I} \right\|_\Ysp \, .
    \label{eq:DecompositionNormAbstract}
\end{equation}

Of course, one has to impose certain restrictions on the covering $\mathcal{Q}$, on the sequence space $\Ysp$ and on the family $(\varphi_i)_{i \in I}$
to ensure that equation (\ref{eq:DecompositionNormAbstract}) yields a well-defined norm/space,
independent of the 
partition of unity $(\varphi_i)_{i \in I}$. 
The most important assumptions are that the covering $\mathcal{Q}$
has the \textbf{finite overlap property} (described below) and that the $(\varphi_i)_{i \in I}$ are
uniformly bounded in the pointwise multiplier algebra of $\Bsp$.

In this paper, however, we will \emph{not} use the general framework of decomposition spaces from \cite{fegr85,fe87}. Instead, we restrict ourselves
to (a slight modification of) the more specialized setting from \cite{boni07}, which we describe now.

We start with an open subset $\mathcal{O}$
of the \emph{frequency space} $\R^d$ and a covering $\mathcal{Q} = (Q_i)_{i \in I}$ of $\mathcal{O}$, which we assume to be of a certain regular form:


\noindent
\begin{defn}\label{defn:IndexClusterAndAdmissibleCovering}(cf.\@ \cite[Def.\@ 2.1 and 2.3]{fegr85}, \cite[Def.\@ 7]{boni07} and \cite[Def.\@ 3.2.8]{VoigtlaenderPhDThesis})
\begin{enumerate}
    \item For $J\subset I$  we define the \textbf{set of neighbors of $J$} as
          \[
          J^{\ast}:=\left\{ i\in I\with\exists j\in J:\, Q_{i}\cap Q_{j}\neq\emptyset\right\} .
          \]
          Inductively, we set
          \[J^{0\ast}:=J \qquad \text{and} \qquad  J^{\left(n+1\right)\ast}:=\left(J^{n\ast}\right)^{\ast} \text{ for } n \in \N_0.\]
          We also define $i^{k\ast}:=\left\{ i\right\} ^{k\ast}$
          for $i\in I$ and $k\in\mathbb{N}_{0}$.
      \item We say that $\mathcal{Q}$ is an \textbf{admissible covering of $\mathcal{O}$} , if $\mathcal{Q}$ is a
            covering of $\mathcal{O}$ with $Q_i \neq \emptyset$ for all $i \in I$ and if the following constant is finite:
            \[
                N_{\mathcal{Q}}:=\sup_{i\in I}\left|i^{\ast}\right| \, .
            \]
            If $N_\calQ < \infty$, we say that $\calQ$ has the \textbf{finite overlap property}.
            \item We say that $\mathcal{Q}$ is an \textbf{almost structured covering} (of $\mathcal{O}$) if it is of the form
                $\left(Q_{i}\right)_{i \in I} = \left(T_{i} Q_i ' +b_{i}\right)_{i\in I}$ for certain $T_i \in {\rm GL}(\R^d)$, $b_i \in \R^d$ and certain
                open sets $Q_i ' \subset \R^d$, satisfying the following properties:
\begin{enumerate}
\item $\calQ$ is an admissible covering of $\calO$,
\item The set $\bigcup_{i \in I} Q_i ' \subset \R^d$ is bounded.
\item The following expression (then a constant) is finite:
\[
C_{\mathcal{Q}}:=\sup_{i\in I}\sup_{j\in i^{\ast}}\left\Vert \smash{T_{i}^{-1}}T_{j}\right\Vert .
\]
\item There is a family $(P_i ')_{i \in I}$ of open sets $P_i ' \subset \R^d$
such that:

    (i)  The sets $\left\{ P_i ' \with i\in I \right\}$ and $\left\{ Q_i ' \with i\in I \right\}$ are both finite.

    (ii) We have $\overline{P_i '} \subset Q_i '$ for all $i \in I$ and
      $\calO = \bigcup_{i \in I} (T_i P_i ' + b_i)$.

\end{enumerate}
\end{enumerate}
\end{defn}

In addition to these assumptions on   $\calQ$ we require
$\Phi = (\varphi_i)_{i \in I}$ to  be uniformly bounded
as pointwise multipliers of $\Bsp$, i.e., we impose
$\| \varphi_i f\|_\Bsp \leq C_\Phi \|f\|_\Bsp$, for all $f \in \Bsp, i \in I.$
For $\Bsp = \mathcal{F}\Lsp^p$,
this boils down to the condition
$ \sup_{i \in I} \| \hat{\varphi_i}\|_\Lisp < \infty.$

The precise definition reads as follows:
\begin{defn}(cf.\@ \cite[Definition 2.2]{fegr85} and \cite[Definition 2]{boni07})
    \label{def:FourierBAPU}

    Let $\mathcal{Q} = (Q_i )_{i \in I}$ be an almost structured covering of $\mathcal{O}$. A family $\Phi = (\varphi_i)_{i \in I}$
    is called a \textbf{bounded admissible partition of unity} (BAPU)
    (\textbf{subordinate to} $\mathcal{Q}$), if the following hold:
    \begin{enumerate}
        \item $\varphi_i \in C_c^\infty (\mathcal{O})$ with $\supp \varphi_i \subset Q_i$ for all $i \in I$,
        \item $\sum_{i \in I} \varphi_i \equiv 1$ on $\mathcal{O}$,
        \item $\sup_{i\in I} \|\mathcal{F}^{-1} \varphi_i \|_{\Lsp^1} <\infty$.
    \end{enumerate}
\end{defn}
\begin{rem*}
In most concrete cases, the requirement $\varphi_i \in C_c^\infty (\calO)$
can be relaxed substantially, as long as
all the involved expressions make sense. In this case, however, the reservoirs $\calD ' (\calO)$ and
$Z'(\calO)$---which will be used in Definitions \ref{def:FourierSideDecompositionSpace} and \ref{def:SpaceSideDecompositionSpaces} below---have
to be replaced by suitable substitutes, cf.\@ \cite[Definition 2.4]{fegr85}.
\end{rem*}

The only remaining assumption which we need to define decomposition spaces pertains
to the sequence space $\Ysp$ from equation (\ref{eq:DecompositionNormAbstract}).
For the sake of simplicity, we restrict ourselves to the case of weighted sequence spaces\footnote{We use the convention $\lsp_u^q (I) = \{(c_i)_{i \in I} \suth (u_i c_i)_{i \in I} \in \lsp^q (I)\}$, with the obvious norm.}, i.e.\@ $\Ysp = \lsp_u^q (I)$ for a certain weight $u = (u_i)_{i \in I}$.
Now, the general theory of decomposition spaces requires $\Ysp$ to be a so-called \textbf{$\calQ$-regular sequence space}, cf.\@ \cite[Definition 2.5]{fegr85}.
In our setting, this leads to the following condition:

\begin{defn}\label{defn:ModerateWeights}
Let $\calQ = (Q_i)_{i \in I}$ be an admissible covering of $\mathcal{O}$. We say that a weight $u = (u_i)_{i \in I}$ is
\textbf{$\calQ$-moderate}, if the constant
\[
C_{u, \calQ} := \sup_{i \in I}\sup_{\ell \in i^\ast } \frac{u_i}{u_\ell}
\]
is finite, i.e.\@ if $u_i \asymp u_\ell$ for $Q_i \cap Q_\ell \neq \emptyset$ (uniformly w.r.t.\@ $i,\ell$).
\end{defn}

Now that we have clarified all of our assumptions, we can finally give a formal definition of the decomposition spaces that we will consider in this paper.
Since the norm on these spaces requires to compute frequency localizations of the form $f_i = \mathcal{F}^{-1} (\varphi_i \widehat{f})$, it is often more convenient to directly \emph{work on the Fourier side}, as in the following definition.
The usual \emph{space sided} version of decomposition spaces will be introduced below.

\begin{defn}
\label{def:FourierSideDecompositionSpace}
Assume that $\Phi=\left(\varphi_{i}\right)_{i\in I}$ is a BAPU subordinate
to the almost structured covering $\mathcal{Q}=\left(Q_{i}\right)_{i\in I}$
of $\mathcal{O}$, let $p,q \in [1,\infty]$, and assume that $u = (u_i)_{i \in I}$ is $\calQ$-moderate.
Then we define for  $f\in\mathcal{D}'\left(\mathcal{O}\right)$
\[
\left\Vert f\right\Vert _{\mathcal{D}_{\mathcal{F}}\left(\mathcal{Q},\Lsp^{p},\lsp_u^q \right)}:=\left\Vert f\right\Vert _{\mathcal{D}_{\mathcal{F},\Phi}\left(\mathcal{Q},\Lsp^{p},\lsp_u^q\right)}:=\left\Vert \left(\left\Vert \mathcal{F}^{-1}\left(\varphi_{i}f\right)\right\Vert _{\Lsp^{p}}\right)_{i\in I}\right\Vert _{\lsp_u^q}\in\left[0,\infty\right],
\]
with the convention that for a family $c=\left(c_{i}\right)_{i\in I}$
with $c_{i}\in\left[0,\infty\right]$, the expression $\left\Vert c\right\Vert _{\lsp_u^q}$
is to be read as $\infty$ if $c_{i}=\infty$ for some $i\in I$ or if $c\notin \lsp_u^q (I)$.

Define the \textbf{Fourier-side decomposition space}
$\mathcal{D}_{\mathcal{F}}\left(\mathcal{Q},\Lsp^{p},\lsp_u^q\right)$ with respect to the
covering $\mathcal{Q}$, integrability exponent $p$ and global component
$\lsp_u^q$ as
\[
    \mathcal{D}_{\mathcal{F}}\left(\mathcal{Q},\Lsp^{p},\lsp_u^q\right):=\left\{ f\in\mathcal{D}'\left(\mathcal{O}\right)\with\left\Vert f\right\Vert _{\mathcal{D}_{\mathcal{F}}\left(\mathcal{Q},\Lsp^{p},\lsp_u^q\right)}<\infty\right\} .
\]
\end{defn}
\begin{rem*}
We observe that $\varphi_{i}f$ is a distribution on $\mathcal{O}$
with compact support in $\calO$, which thus extends to a (tempered) distribution
on $\mathbb{R}^{d}$. By the Paley-Wiener theorem, 
this implies that $\mathcal{F}^{-1}\left(\varphi_{i}f\right)\in\mathcal{S}'\left(\mathbb{R}^{d}\right)$
is given by (integration against) a smooth function. Thus, it makes
sense to write $\left\Vert \mathcal{F}^{-1}\left(\varphi_{i}f\right)\right\Vert _{\Lsp^{p}}$,
with the caveat that this expression could be infinite.

We finally observe that the notations
$\left\Vert \cdot\right\Vert _{\mathcal{D}_{\mathcal{F}}\left(\mathcal{Q},\Lsp^{p},\lsp_u^q\right)}$
and $\mathcal{D}_{\mathcal{F}}\left(\mathcal{Q},\Lsp^{p},\lsp_u^q\right)$ suppress the family
$\left(\varphi_{i}\right)_{i\in I}$ used to define the norm above.
But the general theory of decomposition spaces from \cite{fegr85} implies that different choices yield the same spaces with equivalent norms.

\end{rem*}

Since it is more common to work with the proper, ``space-sided'' objects, rather than with their Fourier-side
versions, we will now define the usual, space-sided version of decomposition spaces. To this end, we first introduce the
reservoir $Z'\left(\mathcal{O}\right)$ which will be used for these spaces. Our notation is inspired by Triebel's book \cite{tr77}.
\begin{defn}
\label{def:UnifiedReservoir}For $\emptyset\neq\mathcal{O}\subset\mathbb{R}^{d}$
open, we define
\[
Z\left(\mathcal{O}\right):=\mathcal{F}\left(C_{c}^{\infty}\left(\mathcal{O}\right)\right):=\left\{ \smash{\widehat{f}}\with f\in C_{c}^{\infty}\left(\mathcal{O}\right)\right\} \leq\mathcal{S}\left(\smash{\mathbb{R}^{d}}\right)
\]
and endow this space with the unique topology that makes the Fourier transform
\[
\mathcal{F}:C_{c}^{\infty}\left(\mathcal{O}\right)\to Z\left(\mathcal{O}\right)
\]
a homeomorphism.

We equip the topological dual space $Z'\left(\mathcal{O}\right):=\left[Z\left(\mathcal{O}\right)\right]'$
of $Z\left(\mathcal{O}\right)$ with the weak-$\ast$-topology, i.e.,
with the topology of pointwise convergence on $Z\left(\mathcal{O}\right)$.

Finally, as on the Schwartz space, we extend the Fourier transform
by duality to $Z'\left(\mathcal{O}\right)$, i.e.\@ we define
\begin{equation}
    \mathcal{F}:Z'\left(\mathcal{O}\right)\to\mathcal{D}'\left(\mathcal{O}\right),f\mapsto \widehat{f} := f\circ\mathcal{F}.\label{eq:FourierTransformOnUnifiedReservoir}
\end{equation}
\end{defn}
\begin{rem*}
Since $\mathcal{F}:C_{c}^{\infty}\left(\mathcal{O}\right)\to Z\left(\mathcal{O}\right)$
is a linear homeomorphism, the Fourier transform as defined in equation
(\ref{eq:FourierTransformOnUnifiedReservoir}) is
a linear homeomorphism as well.
\end{rem*}
Finally, we can define the {\it space-side decomposition spaces}.
\begin{defn}
\label{def:SpaceSideDecompositionSpaces}
In the setting of Definition \ref{def:FourierSideDecompositionSpace} and for $f\in Z'\left(\mathcal{O}\right)$, we define
\[
\left\Vert f\right\Vert_{\mathcal{D}\left(\mathcal{Q},\Lsp^{p},\lsp_u^q\right)} :=\left\Vert \smash{\widehat{f}}\right\Vert _{\mathcal{D}_{\mathcal{F}}\left(\mathcal{Q},\Lsp^{p},\lsp_u^q\right)}=\left\Vert \left(\left\Vert \mathcal{F}^{-1}\left(\varphi_{i}\smash{\widehat{f}}\right)\right\Vert _{\Lsp^{p}}\right)_{i\in I}\right\Vert _{\lsp_u^q}\in\left[0,\infty\right]
\]
and define the \textbf{(space-side) decomposition space}
$\mathcal{D}\left(\mathcal{Q},\Lsp^{p},\lsp_u^q\right)$ with
respect to the covering $\mathcal{Q}$, integrability exponent $p$
and global component $\lsp_u^q$ by
\[
\mathcal{D}\left(\mathcal{Q},\Lsp^{p},\lsp_u^q\right):=\left\{ f\in Z'\left(\mathcal{O}\right)\with\left\Vert f\right\Vert _{\mathcal{D}\left(\mathcal{Q},\Lsp^{p},\lsp_u^q\right)}<\infty\right\} .
\]
\end{defn}
\begin{rem}
\label{rem:TemperedDistributionsAsReservoirIncomplete}Since the Fourier
transform $\mathcal{F}:Z'\left(\mathcal{O}\right)\to\mathcal{D}'\left(\mathcal{O}\right)$
is an isomorphism, it is clear that the Fourier transform restricts
to an (isometric) isomorphism
\[
\mathcal{F}:\mathcal{D}\left(\mathcal{Q},\Lsp^{p},\lsp_u^q\right)\to\mathcal{D}_{\mathcal{F}}\left(\mathcal{Q},\Lsp^{p},\lsp_u^q\right).
\]
Hence, it does not really matter whether one
considers the ``space-side'' or the ``Fourier-side'' version of
these spaces.

At this point, one might ask why we use the spaces $\calD ' (\calO)$ and $Z'(\calO)$ at all, instead of resorting to the more
common reservoir $\SchpRd$, which is for example used to define Besov spaces. The main reasons for this are the
following:

(a)
We want to allow the case $\mathcal{O}\subsetneq\mathbb{R}^{d}$.
If we were to take the space
$\SchpRd$, the decomposition space
norm would \emph{not} be positive definite, or we would have to factor
out a certain subspace of $\SchpRd$.
This is for example done in the usual definition of homogeneous Besov
spaces, which are subspaces of $\SchpRd / \mathcal{P}$,
where $\mathcal{P}$ is the space of polynomials. Here, it seems more
natural to use the space $\mathcal{D}'\left(\mathcal{O}\right)$.

(b)
In case of $\mathcal{O}=\mathbb{R}^{d}$, one could use $\SchpRd$
as the reservoir, as e.g.\@ in \cite{boni07}.
But as shown in \cite[Example 3.4.14]{VoigtlaenderPhDThesis}, this does in general \emph{not} yield a complete space, even for the
case $\Ysp=\lsp_{u}^{1}$ with a $\calQ$-moderate weight $u$.
%
\end{rem}

The following theorem settles the issue of completeness:
\begin{thm}\label{thm:DecompositionSpaceComplete}(cf.\@ \cite[Theorem 1.4.13]{VoigtlaenderPhDThesis})

    Under the assumptions of Definition \ref{def:FourierSideDecompositionSpace}, the (Fourier-side) decomposition space
    $\mathcal{D}_{\mathcal{F}}\left(\mathcal{Q},\Lsp^{p},\lsp_u^q\right)$ is a Banach
    space which embeds continuously into $\mathcal{D}'\left(\mathcal{O}\right)$.

    Likewise, $\mathcal{D}\left(\mathcal{Q},\Lsp^{p},\lsp_u^q\right)$ is also complete and embeds continuously into $Z' (\mathcal{O})$.
\end{thm}

Before we close this section on the basic properties of decomposition spaces, we note that we have always \emph{assumed} that we are given some
BAPU $\Phi = (\varphi_i)_{i \in I}$ subordinate to the covering $\mathcal{Q}$. It is thus important to know whether such a BAPU actually exists.
The next theorem shows that this is the case for every almost structured covering.
We remark that the result itself, and also the proof, are inspired heavily by the construction used in \cite[Proposition 1]{boni07}.

\begin{thm}\label{thm:AlmostStructuredCoveringsAdmitBAPU}(cf.\@ \cite[Theorem 2.8]{vo16})

    Let $\calQ$ be an almost structured covering. Then there is a BAPU $\Phi = (\varphi_i)_{i \in I}$ subordinate to $\calQ$.
\end{thm}

All in all, we now know how to obtain well-defined decomposition spaces with respect to reasonable coverings.
In the next subsection, we consider a special example of decomposition spaces---the $\alpha$ modulation spaces---in greater detail.

\subsection{$\alpha$-modulation spaces}
\label{sub:AlphaModulationSpaces}
The starting point for the original definition of $\alpha$-modulation spaces in Gröbner's thesis\cite{gr92-2} were the two
parallel worlds of Besov-Triebel-Lizorkin spaces and modulation spaces.
Given these two types of spaces, it was natural to ask whether there
is a way to connect these two families in a ``continuous way''.
Although one could of course apply complex interpolation in order
to construct spaces which are (for fixed parameters $p,q,s$) ``in between''  $\BspqRdN$ and $\MspqRdN$,
this approach was---so far---not really successful,
because it seems to be very difficult
to provide constructive characterizations of the members of such spaces
which could be used in practice.

In contrast, it appeared---after some reflection---quite natural to try
to interpolate the two types of spaces \emph{in a geometric sense}, i.e.\@ to consider decompositions of the Fourier
domain which are ``in between'' the \emph{dyadic} partitions of unity, playing
a crucial role in the description of Besov-Triebel-Lizorkin-spaces, and the
\emph{uniform} partitions of unity, which are relevant for the
description of modulation spaces.

The basic observation for this ``geometric interpolation'' approach is that the \emph{uniform} covering $(Q_k)_{k \in \Z^d}$ and the \emph{dyadic} covering $(P_n)_{n \in \N_0}$ satisfy
\begin{align*}
    [\lambda (Q_k)]^{1/d} &\asymp  \langle \xi \rangle^{0} \qquad \forall \xi \in Q_k \qquad \forall k \in \Z^d, \\
    [\lambda(P_n)]^{1/d} &\asymp \langle \xi\rangle^1 \qquad \forall \xi \in P_n \qquad \forall n \in \N_0,
\end{align*}
with $\langle \xi \rangle := 1 + |\xi|$.

Thus, for $\alpha \in [0,1]$, an \textbf{$\alpha$-covering} should be an (open, admissible) covering $(Q_i)_{i \in I}$ of $\R^d$ which satisfies
\[
    [\lambda (Q_i)]^{1/d} \asymp \langle \xi\rangle^{\alpha} \qquad \forall \xi \in Q_i \qquad \forall i \in I,
\]
with the implied constant independent of the precise choice of $\xi$ and $i$.
Apart from this natural condition, a certain further assumption (cf.\@ \cite[Definition 3.1]{boni06-2}) is imposed to rule out pathological coverings.
For brevity, we omit this condition.

Given this definition, one might wonder (at least for $\alpha \in (0,1)$), whether there exist $\alpha$-coverings. This is indeed the case. As shown in
\cite[Proposition A.1]{boni06-2}, there is an $\alpha$ covering $\calQ^{(\alpha)}$ of the form $\calQ^{(\alpha)} = (B_k^r)_{k \in \Z^d \setminus \{0\}}$, with
\[
    B_k^r := \left\{ \xi \in \R^d \with \left| \xi - |k|^{\alpha_0} k\right| < r \cdot |k|^{\alpha_0} \right\},
\]
where $r = r_\alpha >0$ is chosen suitably and where $\alpha_0 := \frac{\alpha}{1-\alpha}$. Furthermore, as shown in \cite[Lemma 6.1.2 and Theorem 6.1.3]{VoigtlaenderPhDThesis},
$\calQ^{(\alpha)}$ is an (almost) structured covering of $\R^d$ and thus admits a BAPU $(\varphi_i)_{i \in I}$.
Finally, in \cite[Lemma B.2]{boni06-2}, it was shown that any two $\alpha$-coverings of $\R^d$ are \emph{equivalent} (in a certain sense).

Given the covering $\calQ^{(\alpha)}$, we need suitable $\calQ^{(\alpha)}$-moderate weights, in order to define the $\alpha$-modulation spaces.
By \cite[Lemma 6.1.2]{VoigtlaenderPhDThesis}, moderateness holds---for arbitrary $\gamma \in \R$---for the weight
\[
    w^{(\gamma)} := (\langle k\rangle^{\gamma})_{k \in \Z^d \setminus \{0\}}\, ,
\]
so that the \textbf{$\alpha$ modulation spaces}
\[
    \Msapq(\Rst^d) := \decompsp{\calQ^{(\alpha)}}{p}{\lsp^q_{w^{(s/(1-\alpha))}}}
\]
are well-defined Banach spaces. Note that the weight $w^{(s/(1-\alpha))}$ satisfies
\[
    \langle \xi \rangle^{s} \asymp w_k^{(s/(1-\alpha))} \qquad \forall k \in \Z^d \setminus \{0\} \, \forall \xi \in B_r^k,
\]
similar to the case of Besov spaces and modulation spaces.  In particular, for the ``limit cases'' $\alpha = 0$ and $\alpha = 1$,
we have $\Msp^{s,0}_{p,q}(\R^d) = \Msp^{s}_{p,q}(\R^d)$, as well as $\Msp^{s, 1}_{p,q}(\R^d) = \Bsp^{s}_{p,q} (\R^d)$.

Since the $\alpha$ modulation spaces are obtained using ``geometric interpolation'' between Besov spaces and modulation spaces, one \emph{could}
expect that $\Msapq(\R^d)$ can also be obtained by ``ordinary'' means of interpolation (like complex interpolation) from Besov- and modulation spaces.
But at least for the case of \emph{complex} interpolation, this is false: In \cite[Theorem 1.1]{faguwuzh15}, it is shown that
\[
    \left[ \Msp^{s_1, \alpha_1}_{p_1, q_1} (\R^d) , \Msp^{s_2, \alpha_2}_{p_2, q_2} (\R^d) \right]_\theta = \Msp^{s, \alpha}_{p,q}(\R^d)
\]
for certain $p, p_1, p_2, q, q_1, q_2 \in [1,\infty]$ and $s, s_1, s_2 \in \R$, as well as $\alpha, \alpha_1, \alpha_2 \in [0,1]$ and $\theta \in (0,1)$
can only hold if
\[
    \alpha_1 = \alpha_2 \quad \text{ or } \quad p_1 = q_1 = 2 \quad \text{ or } \quad p_2 = q_2 = 2.
\]
If one of these conditions hold, then interpolation works as expected. This might be surprising for $\alpha_1 \neq \alpha_2$;
but in this case, we have $p_1 = q_1 = 2$ or $p_2 = q_2 = 2$, which implies $\Msp^{s_1,\alpha_1}_{p_1,q_1} = \Hsp^{s_1} = \Msp^{s_1, \alpha_2}_{p_1, q_1}$ or
$\Msp^{s_2,\alpha_2}_{p_2,q_2} = \Hsp^{s_2} = \Msp^{s_2, \alpha_1}_{p_2, q_2}$, respectively.

But also apart from interpolation, it is natural to ask how the different $\alpha$ modulation spaces are related.
Concretely, one might wonder under which conditions an embedding of the form
\begin{equation}
    \Msp^{s_1, \alpha_1}_{p_1, q_1} (\R^d) \hookrightarrow \Msp^{s_2, \alpha_2}_{p_2, q_2}(\R^d)
    \label{eq:AlphaModulationEmbedding}
\end{equation}
holds. For $\alpha_1, \alpha_2 \in \{0,1\}$, this question was solved by Kobayashi and Sugimoto\cite{kosu11}.
For general $\alpha$, it was considered by Toft and Wahlberg\cite{towa12}, shortly before it was solved completely---for $(p_1, q_1) = (p_2, q_2)$---by
Han and Wang\cite{hawa14}. Based in part on their ideas, the second named author
of the present paper developed a general theory of embeddings between decomposition spaces (cf.\@ \cite{VoigtlaenderPhDThesis}), which we will outline in
Subsection \ref{sub:DecompositionEmbedding}. 
Using this theory, the above question can be solved easily---even for $(p_1, q_1) \neq (p_2, q_2)$---cf.\@ \cite[Theorems 6.1.7 and 6.2.8]{VoigtlaenderPhDThesis}:
\begin{thm}
    \label{thm:AlphaModulationEmbeddings}
    Let $0 \leq \alpha \leq \beta \leq 1$, $p_1, p_2, q_1, q_2 \in [1,\infty]$ and $s_1, s_2 \in \R$. Define
    \begin{align*}
        s^{(0)} := \alpha \left( \frac{1}{p_2} - \frac{1}{p_1} \right) + (\alpha - \beta) \left( \frac{1}{p_2^{\triangledown}} - \frac{1}{q_1} \right)_+ \quad ,\\
        s^{(1)} := \alpha \left( \frac{1}{p_2} - \frac{1}{p_1} \right) + (\alpha - \beta) \left( \frac{1}{q_2} - \frac{1}{p_1^\triangle} \right)_+ \quad,
    \end{align*}
    with $p^{\triangledown} := \min \{p, p'\}$, as well as $p^{\triangle} := \max \{p, p'\}$.

    We have
    \[
        \Msp^{s_1, \alpha}_{p_1, q_1} (\R^d) \hookrightarrow \Msp^{s_2, \beta}_{p_2, q_2}(\R^d)
    \]
    if and only if $p_1 \leq q_1$ and
    \[
        \begin{cases}
            s_2 \leq s_1 + d \cdot s^{(0)} ,                                                                       & \text{if } q_1 \leq q_2,\\
            s_2 <    s_1 + d \cdot \left( s^{(0)} + (1-\beta)\left( q_1^{-1} - q_2^{-1} \right) \right), & \text{if } q_1 > q_2.
        \end{cases}
    \]

    Conversely, we have
    \[
        \Msp^{s_1, \beta}_{p_1, q_1} (\R^d) \hookrightarrow \Msp^{s_2, \alpha}_{p_2, q_2}(\R^d)
    \]
    if and only if $p_1 \leq q_1$ and
    \[
        \begin{cases}
            s_2 \leq s_1 + d \cdot s^{(1)} ,                                                                           & \text{if } q_1 \leq q_2, \\
            s_2 <    s_1 + d \cdot \left( s^{(1)} + (1 - \beta) \left(  q_1^{-1} - q_2^{-1} \right) \right), & \text{if } q_1 > q_2.
        \end{cases}
    \]
\end{thm}
Using the theory of embeddings between decomposition spaces from Subsection \ref{sub:DecompositionEmbedding}, one can explain the main geometric properties
of the $\alpha$ coverings $\calQ^{(\alpha)}$ which lead to the preceding theorem: The main point is that the covering $\calQ^{(\alpha)}$ is
\emph{almost subordinate to} (``finer than'', cf.\@ equation (\ref{eq:QAlmostSubordinateToP})) the covering $\calQ^{(\beta)}$ for $\alpha \leq \beta$.
Furthermore, the precise conditions depend on the number of ``smaller sets'' that are needed to cover the ``bigger'' sets,
cf.\@ \cite[Lemma 6.1.5]{VoigtlaenderPhDThesis} and Theorems \ref{thm:DecompositionEmbeddingPIntoQ} and \ref{thm:DecompositionEmbeddingQIntoP} from below.
We finally remark that the theorems in \cite{VoigtlaenderPhDThesis} apply for the full range $(0,\infty]$ of the exponents.
But in the present paper, we restrict ourselves to the range $[1,\infty]$ for simplicity.

\subsection{Embeddings between different decomposition spaces}
\label{sub:DecompositionEmbedding}

In this subsection, we consider embeddings between decomposition spaces with respect to different coverings, i.e.\@ of the form
\begin{equation}
    \decompsp{\calQ}{p_1}{\lsp_u^{q_1}} \hookrightarrow \decompsp{\calP}{p_2}{\lsp_v^{q_2}}
    \label{eq:DecompositionEmbeddingQIntoP}
\end{equation}
for $p_1, p_2, q_1, q_2 \in [1,\infty]$ and two almost structured coverings
\[\calQ = (Q_i)_{i \in I} = (T_i Q_i ' + b_i)_{i \in I} \qquad \text{ and } \qquad \calP = (P_j)_{j \in J} = (S_j P_j ' + c_j)_{j \in J}\]
of two (possibly different) subsets $\calO , \calO '$ of the frequency space $\R^d$.
We assume the weights $u = (u_i)_{i \in I}$ and $v = (v_j)_{j \in J}$ to
be moderate w.r.t.\@ $\calQ$ and $\calP$, respectively.

As seen in the previous subsection, our setting includes embeddings between $\alpha$ modulation spaces for different values of $\alpha$.

As the main standing requirement for this subsection, we assume that $\calQ$ is \textbf{almost subordinate} to $\calP$.
Very roughly, this means that the sets $Q_i$ are ``smaller'' than the sets $P_j$, or that $\calQ$ is ``finer'' then $\calP$.
Rigorously, it means that there is some fixed $n \in \N_0$
such that for every $i \in I$, there is some $j_i \in J$ satisfying
\begin{equation}
    Q_i  \subset P_{j_i}^{n\ast} := \bigcup_{\ell \in j_i^{n\ast}} P_\ell.
    \label{eq:QAlmostSubordinateToP}
\end{equation}
Note that this assumption implies $\calO \subset \calO '$. Even more importantly, this assumption destroys the symmetry between $\calQ, \calP$
in equation (\ref{eq:DecompositionEmbeddingQIntoP}), so that in addition to (\ref{eq:DecompositionEmbeddingQIntoP}), we will also consider
the ``reverse'' embedding
\begin{equation}
    \decompsp{\calP}{p_2}{\lsp_v^{q_2}} \hookrightarrow   \decompsp{\calQ}{p_1}{\lsp_u^{q_1}}.
    \label{eq:DecompositionEmbeddingPIntoQ}
\end{equation}

Under the assumption that $\calQ$ is almost subordinate to $\calP$, the object which describes those features of
the coverings $\calQ, \calP$ which are relevant for us is the family of \textbf{intersection sets} given by
\[
    I_j := \left\{ i \in I \with Q_i \cap P_j \neq \emptyset \right\} \qquad \text{ for } j \in J.
\]
Of course, $I_j = \emptyset$ if and only if $P_j \cap \calO = \emptyset$. Since these sets will be irrelevant for our purposes, we additionally define
\[
    J_{\calO} := \left\{ j \in J \with P_j \cap \calO \neq \emptyset \right\}.
\]

In the remainder of this subsection, we will state sufficient conditions and necessary conditions for the existence of the embeddings
(\ref{eq:DecompositionEmbeddingQIntoP}) and (\ref{eq:DecompositionEmbeddingPIntoQ}), respectively.
In general, these two conditions will only coincide for a certain range of $p_1$ or $p_2$, while
there is a gap between the two conditions outside of this range.

Under suitable additional assumptions, however, more strict necessary conditions can be derived;
in fact, a \emph{complete characterization} of the existence of the embeddings can be achieved. For this to hold, we will
(occasionally, but not always) assume that the following properties are fulfilled:
\begin{defn}\label{def:RelativeModerateness}
    \begin{enumerate}
        \item We say that the weight $u = (u_i)_{i \in I}$ is \textbf{relatively $\calP$-moderate}, if
            \[
                \sup_{j \in J} \,\, \sup_{i,\ell \in I_j} \frac{u_i}{u_\ell} < \infty.
            \]
        \item The (almost structured) covering $\calQ = (T_i Q_i ' + b_i)_{i \in I}$ is called \textbf{relatively $\calP$-moderate}, if
            the weight $(|\det T_i|)_{i \in I}$ is relatively $\calP$-moderate.

            Roughly speaking, this means that two (small) sets $Q_i, Q_\ell$ have essentially the same measure if they intersect the same (large)
            set $P_j$.

    \end{enumerate}
\end{defn}
Although these assumptions might appear rather restrictive, they are fulfilled in many practical cases; in particular if $\calQ$ and $\calP$ are coverings
associated to $\alpha$-modulation spaces, and if $u,v$ are the usual weights for these spaces.

We can now analyze existence of the embedding (\ref{eq:DecompositionEmbeddingPIntoQ}):
\begin{thm}(cf.\@ \cite[Theorem 5.4.1]{VoigtlaenderPhDThesis})
    \label{thm:DecompositionEmbeddingPIntoQ}

    Define
        $p_2^\triangle := \max \{p_2, p_2 '\}$
    and for $r \in [1,\infty]$, let
    \begin{align*}
        C_1^{(r)} &:= \left\| \left( \left\| \left( |\det T_i|^{\frac{1}{p_2} - \frac{1}{p_1}} \cdot u_i \right)_{i \in I_j}\right\|_{\lsp^{q_1 \cdot (r / q_1)'} (I_j)} \bigg/  v_j \right)_{\!\!j \in J}\right\|_{\lsp^{q_1 \cdot (q_2 / q_1)'}(J)} , \\
        C_2 &:= \left\| \left( \frac{u_{i_j}}{v_j} \cdot  |\det T_{i_j}|^{\frac{1}{p_2} - \left( \frac{1}{q_1} - \frac{1}{p_2^\triangle} \right)_{\!\!+} - \frac{1}{p_1} } \cdot  |\det S_j|^{\left( \frac{1}{q_1} - \frac{1}{p_2^\triangle} \right)_{\!\!+}}    \right)_{\!\!j \in J_{\calO}}\right\|_{\lsp^{q_1 \cdot (q_2 / q_1)'}(J_{\calO})} ,
    \end{align*}
    where for each $j \in J_{\calO}$ some $i_j \in I$ with $Q_{i_j} \cap P_j \neq \emptyset$ can be chosen arbitrarily.

    Then the following hold:
    \begin{enumerate}
        \item If $C_1^{(p_2^\triangle)} < \infty$ and $p_2 \leq p_1$, then the map
              \begin{equation}
                  \iota : \decompsp{\calP}{p_2}{\lsp_v^{q_2}} \to \decompsp{\calQ}{p_1}{\lsp_u^{q_1}}, f \mapsto f|_{\mathcal{F}(C_c^\infty (\calO))}
                  \label{eq:ExplicitEmbeddingPIntoQ}
              \end{equation}
              is well-defined and bounded.

        \item Conversely, if
            \begin{equation}
                \qquad \qquad \theta : \left(\mathcal{F}^{-1}(C_c^\infty (\calO)) , \fsn {\cdot} {\decompsp{\calP}{p_2}{\lsp_v^{q_2}}}\right) \to \decompsp{\calQ}{p_1}{\lsp_u^{q_1}}, f \mapsto f
                \label{eq:WeakEmbeddingPIntoQ}
            \end{equation}
            is bounded, then $C_1^{(p_2)} < \infty$ and $p_2 \leq p_1$.

        \item Finally, if $P_j \subset \calO$ holds\footnote{The main case in which this holds is if $\calO = \calO'$.} for all $j \in J_{\calO}$
              and if additionally the covering $\calQ$ and the weight $u = (u_i)_{i \in I}$ are both relatively $\calP$-moderate,
              we have the following equivalence:
              \begin{align*}
                  \iota \text{ is bounded} & \quad \Longleftrightarrow \quad \theta \text{ is bounded} \\
                                           & \quad \Longleftrightarrow \quad \left(C_1^{(p_2^\triangle)} < \infty \text{ and } p_2 \leq p_1 \right) \\
                                           & \quad \Longleftrightarrow \quad \left(C_2 < \infty \text{ and } p_2 \leq p_1 \right).
              \end{align*}
    \end{enumerate}
\end{thm}
\begin{rem*}
    \begin{enumerate}
        \item We achieve a \emph{complete characterization} of the existence of the embedding (\ref{eq:DecompositionEmbeddingPIntoQ}) if $\calQ$ and $u$ are relatively $\calP$-moderate, but also in case of $p_2 \in [2,\infty]$, since in this case, $p_2^{\triangle} = p_2$ and hence \@ $C_1^{(p_2^\triangle)} = C_1^{(p_2)}$.

        \item Even for well-understood special cases like $\alpha$ modulation spaces, the above theorem yields new results, since even in the most general previous result \cite{hawa14}, only the case $(p_1, q_1) = (p_2, q_2)$ was studied.

        \item Considering $\iota$ as an \emph{embedding} is not always justified. For example, if $\calO' \setminus \calO$ has nonempty interior,
              then every $f \in \mathcal{F}^{-1} ( C_c^\infty (\calO' \setminus \calO))$ satisfies $\iota f = 0$, although $f \neq 0$ is possible.
    \end{enumerate}
\end{rem*}

The result for the embedding (\ref{eq:DecompositionEmbeddingQIntoP}) is similar:
\begin{thm}(cf.\@ \cite[Theorem 5.4.4]{VoigtlaenderPhDThesis})
    \label{thm:DecompositionEmbeddingQIntoP}

    Let $(\varphi_i)_{i \in I}$ be a BAPU for $\calQ$. Define
        $p_2^\triangledown := \min \{p_2, p_2 '\}$
    and for $r \in [1,\infty]$, let
    \begin{align*}
        C_1^{(r)} &:= \left\| \left( v_j \cdot \left\| \left(u_i^{-1} \cdot  |\det T_i|^{\frac{1}{p_1} - \frac{1}{p_2}} \right)_{i \in I_j}\right\|_{\lsp^{r \cdot (q_1 / r)'} (I_j)} \right)_{\!\!j \in J}\right\|_{\lsp^{q_2 \cdot (q_1 / q_2)'}(J)} , \\
        C_2 &:= \left\| \left(  \frac{v_j}{u_{i_j}} \cdot |\det T_{i_j}|^{\frac{1}{p_1} - \left( \frac{1}{p_2^{\triangledown}} - \frac{1}{q_1} \right)_{\!\! +} - \frac{1}{p_2}} \cdot  |\det S_j|^{\left( \frac{1}{p_2^{\triangledown}} - \frac{1}{q_1} \right)_{\!\! +}} \right)_{\!\!j \in J_{\calO}}\right\|_{\lsp^{q_2 \cdot (q_1 / q_2)'}(J_{\calO})} ,
    \end{align*}
    where for each $j \in J_{\calO}$ some $i_j \in I$ with $Q_{i_j} \cap P_j \neq \emptyset$ can be chosen arbitrarily.

    Then the following hold:
    \begin{enumerate}
        \item If $C_1^{(p_2^\triangledown)} < \infty$ and if $p_1 \leq p_2$, then the map
            \begin{equation}
                \quad \qquad \iota = \iota_{\Phi} : \decompsp{\calQ}{p_1}{\lsp_u^{q_1}} \hookrightarrow \decompsp{\calP}{p_2}{\lsp_v^{q_2}}, f \mapsto \sum_{i \in I} \mathcal{F}^{-1} ( \varphi_i \cdot \widehat{f}\,)
                \label{eq:ExplicitEmbeddingQIntoP}
            \end{equation}
            is well-defined and bounded. 
            Here, $\iota f$ acts as follows:
            \[
                \langle \iota f, \gamma \rangle = \sum_{i \in I} \langle f, \mathcal{F} (\varphi_i \cdot \mathcal{F}^{-1} \gamma) \rangle \quad \text{ for } \quad \gamma \in \mathcal{F}(C_c^\infty (\calO ')),
            \]
            with absolute convergence of the series for $f \in \decompsp{\calQ}{p_1}{\lsp_u^{q_1}}$.

            Furthermore, $\langle \iota f, \gamma \rangle = \langle f, \gamma \rangle$ holds for all $\gamma \in \mathcal{F}(C_c^\infty (\calO))$, so that
            $\iota f \in \decompsp{\calP}{p_2}{\lsp_v^{q_2}} \subset Z'(\calO ')$ extends $f \in \decompsp{\calQ}{p_1}{\lsp_u^{q_1}} \subset Z'(\calO)$.


        \item If the map
            \begin{equation}
                \quad \qquad \theta : \left( \mathcal{F}^{-1}(C_c^\infty (\calO)), \fsn {\cdot} {\decompsp{\calQ}{p_1}{\lsp_u^{q_1}}}\right) \to \decompsp{\calP}{p_2}{\lsp_v^{q_2}}, f \mapsto f
                \label{eq:WeakEmbeddingQIntoP}
            \end{equation}
            is bounded, then $p_1 \leq p_2$ and $C_1^{(p_2)} < \infty$.

        \item Finally, if $P_j \subset \calO$ holds for all $j \in J_{\calO}$ and if additionally the covering $\calQ$ and the weight $u = (u_i)_{i \in I}$
              are both relatively $\calP$-moderate, we have the following equivalence:
              \begin{align*}
                  \iota \text{ is bounded} & \quad \Longleftrightarrow \quad \theta \text{ is bounded}. \\
                                           & \quad \Longleftrightarrow \quad \left(C_1^{(p_2^\triangledown)} < \infty \text{ and } p_1 \leq p_2 \right) \\
                                           & \quad \Longleftrightarrow \quad \left(C_2 < \infty \text{ and } p_1 \leq p_2 \right).
              \end{align*}
    \end{enumerate}
\end{thm}
\begin{rem*}
    \begin{enumerate}
        \item We achieve a complete characterization of the existence of the embedding (\ref{eq:DecompositionEmbeddingQIntoP}) for $p_2 \in [1,2]$. If $\calQ$ and $u$ are relatively $\calP$-moderate, we get a complete characterization for arbitrary $p_2 \in [1,\infty]$.

        \item In contrast to Theorem \ref{thm:DecompositionEmbeddingPIntoQ}, $\iota$ is always injective in the present setting.

        \item Note that the definition of $\iota$ is independent of $p_1, p_2, q_1, q_2$ and $u,v$. In fact, if $\calO = \calO '$, then $\iota f = f$ for all $f \in \decompsp{\calQ}{p_1}{\lsp_u^{q_1}}  \subset Z'(\calO) = Z'(\calO ')$.
    \end{enumerate}
\end{rem*}

For one concrete application of Theorems \ref{thm:DecompositionEmbeddingPIntoQ} and \ref{thm:DecompositionEmbeddingQIntoP}, we refer the reader to the
characterization of embeddings between different $\alpha$ modulation spaces in Theorem \ref{thm:AlphaModulationEmbeddings}. Further applications will be given
in Section \ref{sec:WaveletTheory}.

There are also results which apply if neither $\calQ$ is almost subordinate to $\calP$, nor vice versa.
In fact, it suffices if one can write $\calO \cap \calO ' = A \cup B$, such that $\calQ$ is almost subordinate to $\calP$ \emph{near $A$}
and vice versa \emph{near $B$}. For a precise formulation of this condition, and the resulting embedding results,
we refer to \cite[Theorem 5.4.5]{VoigtlaenderPhDThesis}.

Finally, there are also results for embeddings of the form
\[
    \decompsp{\calQ}{p}{\Ysp} \hookrightarrow W^{k,q}(\R^d)
\]
for the classical Sobolev spaces $W^{k,q}$. These results are easy to apply, since no subordinateness is required.
As shown in \cite{vo16}, existence of the embedding can be completely characterized for $q \in [1,2] \cup \{\infty\}$, while for $q \in (2,\infty)$,
certain sufficient and certain necessary criteria are given; but in general, these do not coincide.

\subsection{Banach frames for Decomposition spaces}
\label{sub:BanachFramesForDecompositionSpaces}
Borup and Nielsen\cite{boni07} gave a construction of Banach frames for decomposition spaces which applies in a very general setting.
In particular, their  
 construction applies to ($\alpha$)-modulation spaces and Besov spaces.
Since this construction makes the power of Banach frames available
for decomposition spaces, we could not resist discussing their results.

Furthermore, their results fit well into the present context: As Borup and Nielsen write themselves: ``[our] frame expansion
should perhaps be considered an adaptable variant of the $\varphi$-transform of Frazier and Jawerth'' (cf.\@ \cite[Section 3.2]{boni07}).


To describe their construction, we first introduce 
\textbf{structured} coverings. An almost structured
covering $\calQ = (Q_i)_{i \in I} = (T_i Q_i ' + b_i)_{i \in I}$ of $\calO = \R^d$ is called structured if $Q_i ' = Q$ for all $ i\in I$, i.e., if all
$Q_i$ are affine images of a fixed set.

The idea 
is to transfer the orthonormal basis $(e^{2 \pi i \langle k , \cdot \rangle })_{k \in \Z^d}$
of $L^2 ([-1/2 ,1/2]^d)$ to each of the sets $Q_i^{(a)} = T_i [-a,a]^d + b_i \supset Q_i$ for certain $a > 0$.
Then, one truncates these periodic functions using a certain (quadratic)
partition of unity subordinate to $\calQ$. Thus, one obtains a tight frame for $\Lsp^2 (\R^d)$. The nontrivial part is to show
that one also obtains Banach frames for the whole range of decomposition spaces.

The construction proceeds as follows:  By \cite[Proposition 1]{boni07}
on finds a family $(\theta_i)_{i \in I}$ such that:
\begin{enumerate}
    \item $\supp \theta_i \subset Q_i$ for all $i \in I$,
    \item $\sum_{i \in I} \theta_i^2 \equiv 1$ on $\calO = \R^d$,
    \item $\sup_{i \in I} \|\mathcal{F}^{-1} \theta_i \|_{\Lsp^1} < \infty$,
    \item $\sup_{i \in I} \|\partial^{\alpha} [\theta_i (T_i \cdot + b_i)] \|_{\sup} <\infty$ for all $\alpha \in \N_0^d$.
\end{enumerate}

Given such a family $(\theta_i)_{i \in I}$, we choose a cube $Q_a \subset \R^d$ of side-length $2a$ satisfying $Q \subset Q_a$.
Finally, for $i \in I$ and $n \in \Z^d$, define $e_{n,i} : \R^d \to \Compl$ by
\[
    e_{n,i} (\xi) := (2a)^{-d/2} \cdot |\det T_i|^{-1/2} \cdot \chi_{Q_a} (T_i^{-1} (\xi - b_i)) \cdot e^{i \frac{\pi}{a} n \cdot T_i^{-1}(\xi - b_i)} \text{ for } \xi \in \R^d,
\]
and set
\[
    \eta_{n,i} := \mathcal{F}^{-1} (\theta_i \cdot e_{n,i}).
\]

Since the family $(e_{n,i})_{n \in \Z^d}$ forms an orthonormal basis of $\Lsp^2(T_i Q_a + b_i)$ and because of $\sum_{i \in I} \theta_i^2 \equiv 1$, it follows
(cf.\@ \cite[Proposition 2]{boni07}) that the family $(\eta_{n,i})_{n\in \Z^d, i\in I}$ forms a \emph{tight frame} for $\Lsp^2 (\R^d)$.

Of course, we are not simply interested in (tight) frames for $\Lsp^2 (\R^d)$ with a given form of frequency localization---we want to obtain a Banach frame
for the decomposition space $\decompsp{\calQ}{p}{\lsp_u^q}$. Thus, we define the $\Lsp^p$-normalized version
\[
    \eta_{n,i}^{(p)} := |\det T_i|^{\frac{1}{2} - \frac{1}{p}} \cdot \eta_{n,i} \quad \text{ for } i \in I \text{ and } n \in \Z^d.
\]
Then, Borup and Nielsen showed (cf.\@ \cite[Proposition 3, Definition 8, Lemma 4 and Theorem 2]{boni07})
that there is a suitable solid BK space $d(\calQ, \lsp^p, \lsp_u^q)$ such that the \emph{coefficient operator}
\[
    \coeff : \decompsp{\calQ}{p}{\lsp_u^q} \to d(\calQ, \lsp^p, \lsp_u^q), f \mapsto (\langle f, \eta_{n,i}^{(p)}\rangle)_{n \in \Z^d, i\in I}
\]
is bounded. As familiar by now, there is also a bounded \emph{reconstruction operator} $\recon : d(\calQ, \lsp^p, \lsp_u^q) \to \decompsp{\calQ}{p}{\lsp_u^q}$
which satisfies $\recon \circ \coeff = {\rm{id}}_{\decompsp{\calQ}{p}{\lsp_u^q}}$. Thus, the family $(\eta_{n,i}^{(p)})_{i \in I, n\in \Z^d}$ forms a
Banach frame for $\decompsp{\calQ}{p}{\lsp_u^q}$.
\section{Abstract and Concrete Wavelet Theory}
\label{sec:WaveletTheory}

As noted in Section \ref{sec:CoorbitTheory}, the description of the spaces $\Bspq$ and $\Fspq$
via the $\varphi$-transform---or via wavelets---can be viewed (at least in part)
as an application of the general theory of coorbit spaces to the affine group which acts on $\Lsp^2(\R^d)$ via translations and \emph{isotropic} dilations. To be precise, this group action yields the \emph{homogeneous} Besov- and Triebel-Lizorkin spaces $\BspqHom$ and $\FspqHom$, not the inhomogeneous ones.

One well known characterization of (homogeneous) Besov spaces shows that these spaces are obtained by placing certain integrability conditions
on the \textbf{(continuous) wavelet transform}
\[
    W_{\varphi} f : \R^d \times \R^{\ast} \to \Compl, (x,a) \mapsto \langle f, T_x D_a \varphi \rangle
\]
of a function or distribution $f$ and a certain \textbf{analyzing window} $\varphi$.
Other applications of the wavelet transform include the characterization of the \textbf{wave-front} set of a distribution using the decay
of the transform \cite{pivu06}
However, due to the isotropic nature of the dilations
\[
    D_a f(x) = a^{-d/2} \cdot f(a^{-1}x),
\]
such a \emph{single wavelet} characterization is only valid in dimension $d=1$ (cf.\@ \cite{cado05} or \cite[Lemma 4.4, Lemma 4.10]{fefuvo14}),
where smoothness is an ``undirected property''.

Even beyond this specific problem of characterizing the wave-front set, it was noted in recent years that the isotropic, directionless nature of the wavelet
transform is a limitation for many applications. To overcome this problem, a vast number of ``directional'' variants of wavelets were invented: In particular,
curvelets\cite{cado05,cado05-1} and shearlets\cite{kula09,dakustte09}.
Among these two systems, shearlets have the special property that there is---as in the case of wavelets---an underlying dilation
group through which the family of shearlets can be generated from a single ``mother wavelet'', see also Section \ref{sec:BeyCoorbitTheory}.

In view of these two very different dilation groups---the affine group and the shearlet group---it becomes natural to consider the bigger picture: Given any
(closed) subgroup $H \subset \GL(\R^d)$, one can form the group $G = \R^d \rtimes H$ of all affine mappings generated by arbitrary translations
and all dilations in $H$. The multiplication on $G$ is given by $(x,h)(y,g) = (x+hy, hg)$ and the Haar measure is
\begin{equation}
    \diff(x,h) = |\det h|^{-1} \diff x \, \diff h,
    \label{eq:HaarMeasureSemidirectProduct}
\end{equation}
where $\diff h$ is the Haar measure of $H$.

The group $G$ from above acts unitarily on $\Lsp^2 (\R^d)$ via translations and dilations, i.e., by the \textbf{quasi-regular representation}
\begin{equation}
    \pi : G \to \calU(\Lsp^2(\R^d)), (x,h) \mapsto T_x D_h.
    \label{eq:QuasiRegularRepresentationDefinition}
\end{equation}
This representation comes with an associated
\textbf{(generalized) wavelet transform}
\[
    W_\varphi f : G \to \Compl, (x,h) \mapsto \langle f, \pi(x,h) \varphi \rangle \quad \text{ for } \quad f,\varphi \in \Lsp^2 (\R^d),
\]
where the (fixed) function $\varphi$ is called the \textbf{analyzing window}. In the general description of coorbit theory in
Section \ref{sec:CoorbitTheory}, this was called the \emph{voice transform}.

Given this transform, it is natural to ask 
which properties of $f$ can be easily read off from $W_\varphi f$.
As for wavelets (for $d=1$ \cite{pivu06})
 or for shearlets (for $d=2$ \cite{kula09,gr11-8}),
 it turns out\cite{fefuvo14} that for
large classes of dilation groups, the wave-front set of a (tempered) distribution $f$ can be characterized via the decay of
$W_\varphi f$.

Another important property of a generalized wavelet system (like shearlets) are its \emph{approximation theoretic properties}.
Here, the question is: Which classes of functions can be approximated well by linear combinations of only a few elements of the wavelet system?
For ``ordinary'' wavelets, this question leads to the theory of Besov spaces and their atomic decompositions, as explored by Frazier and Jawerth.
For a general dilation group, these approximation theoretic properties are (at least in part) encoded by the associated \textbf{wavelet type coorbit spaces},
which we will now discuss in greater detail. One particular problem which is of interest to us is the following:
\emph{If a function/signal $f$ can be well approximated by one wavelet system, can it also be well approximated using a different wavelet system?} Of course,
the answer to this question will depend on the precise nature of the two wavelet systems and on the way in which the statement
``$f$ can be well approximated by \ldots'' is made mathematically precise.

\subsection{General wavelet type coorbit spaces}
\label{sub:WaveletCoorbitTheory}
As long as $\pi$ acts \emph{irreducibly} on $\Lsp^2$ (and if the representation is (square) integrable),
we can apply the general coorbit theory as described in Section \ref{sec:CoorbitTheory} to form the coorbit spaces
\begin{equation}
    \Coosp (G, \Ysp) = \left\{ f \in \calR \with W_\varphi f \in \Ysp \right\},
    \label{eq:CoorbitSpaceAbstractDefinition}
\end{equation}
where $\calR = \calR_\Ysp$ is a suitable \textbf{reservoir},
which plays the role of $\ScRdp$ for the usual Besov- or Triebel-Lizorkin spaces.

Furthermore, $\varphi \in \Lsp^2(\R^d)$ has to be a suitable \textbf{analyzing window}.
Formally, this means that $\varphi$ must fulfill
$\varphi \in \calA_{v_0}$ (cf.\@ equation (\ref{eq:AnalyzingWindowClass})) for a so-called \textbf{control weight} $v_0 : G \to (0,\infty)$.
As we will see (cf.\@ Definition \ref{defn:VanishingMoments}), this condition is closely related to the
usual ``vanishing moments condition'' for ordinary wavelets.

In this section, instead of the general coorbit spaces $\Coosp (\Ysp) = \Coosp(G, \Ysp)$, we will consider the more restrictive case
of the (weighted) mixed Lebesgue space $\Ysp = \Lsp_m^{p,q}(G)$ for $p,q \in [1,\infty]$ and a weight $m : H \to (0,\infty)$.
Precisely, the space $\Lsp_m^{p,q}(G)$ is the space of all measurable functions $f : G \to \Compl$ for which the norm
\[
    \|f\|_{\Lsp_m^{p,q}} := \left\| h \mapsto m(h) \cdot \|f(x,h)\|_{\Lsp^p (\R^d, \d x)} \right\|_{\Lsp^q (H, \, \diff h / |\det h|)}
\]
is finite. This normalization---in particular the measure $\diff h / |\det h|$ on $H$---is chosen such that we have
$\Lsp^{p,p}_m (G) = \Lsp_m^p (G)$, cf.\@ equation (\ref{eq:HaarMeasureSemidirectProduct}).

Recall from Section \ref{sec:CoorbitTheory} that the space $\Ysp$ needs to be a solid BF space on $G$ which is invariant under left- and right translations.
Clearly, $\Ysp = \Lsp_m^{p,q}(G)$ satisfies all of these properties, except possibly for invariance under left- and right translations. To ensure this, we assume
that $m$ is \textbf{$v$-moderate} for some (measurable, locally bounded, submultiplicative) weight $v : H \to (0,\infty)$, i.e., we assume
\[
    m(xyz) \leq v(x) m(y) v(z) \qquad \forall x,y,z \in H.
\]
Under these assumptions, it is shown in \cite[Lemma 1 and Lemma 4]{fuvo15} that there is a \textbf{control weight}
$v_0 : G \to (0,\infty)$ for $\Ysp = \Lsp_m^{p,q}(G)$ which is (by abuse of notation) of the form
$v_0 (x,h) = v_0 (h)$ and measurable, submultiplicative and locally bounded.

In the present setting, coorbit theory can be seen as a theory of \emph{nice wavelets} and \emph{nice signals}, cf.\@ \cite{FuehrWaveletCoorbitOverviewSampta}.
Nice wavelets are those belonging to the class $\calA_{v_0}$ of analyzing windows, while a nice signal $f$ (with respect to an analyzing wavelet $\varphi$)
is one for which $W_\varphi f \in \Lsp_m^{p,q}$, i.e., for which $f \in \Coosp (\Lsp_m^{p,q})$. Now, coorbit theory---if it is applicable---yields
two main properties:
\begin{itemize}
    \item A \emph{consistency statement}: Nice wavelets agree upon nice signals,
        i.e.\@ if $\varphi, \psi \in \calA_{v_0} \setminus \{0\}$ (with $v_0$ depending on $p,q,m$), then
          \[
              W_\varphi f \in \Lsp_m^{p,q} \qquad \Longleftrightarrow \qquad W_\psi f \in \Lsp_m^{p,q}.
          \]
          We even get a norm-equivalence, so that the coorbit space from equation (\ref{eq:CoorbitSpaceAbstractDefinition}) with $\Ysp = \Lsp_m^{p,q}(G)$
          is well-defined.

      \item An \emph{atomic decomposition result}: As seen in Theorem \ref{thm:CoorbitAtomicDecomposition}, we can guarantee atomic
            decompositions of the form
                $f = \sum_{g \in G_0} \alpha_g (f) \cdot \pi(g) \psi$
            for elements $f \in \Coosp(\Lsp_m^{p,q})$ and coefficients $(\alpha_g (f))_{g \in G_0}$ lying in a suitable sequence space,
            if $\psi$ is a \textbf{better vector} (in comparison to just being an analyzing vector), i.e.\@ if $\psi \in \calB_{v_0}$.
%
\end{itemize}

Despite these pleasant features, the theory of (generalized) wavelet type coorbit spaces raises several questions:
\begin{enumerate}[label=(Q\arabic*),ref=(Q\arabic*)]
    \item\label{item:DiscreteSeriesQuestion} For which dilation groups $H$ is the quasi-regular representation from equation (\ref{eq:QuasiRegularRepresentationDefinition})
        irreducible and (square)-integrable, so that coorbit theory is applicable \emph{in principle}?

    \item \label{item:AnalyzingVectorsQuestion} Is coorbit theory applicable, i.e.\@ are there ``nice wavelets''?
          Precisely, do we have $\calA_{v_0} \neq \{0\}$ and $\calB_{v_0} \neq \{0\}$ and are there convenient sufficient criteria for a function
          $\varphi \in \Lsp^2 (\R^d)$ to belong to $\calA_{v_0}$ or to $\calB_{v_0}$?

      \item \label{item:CoorbitSpacesEmbeddingQuestion} How are the resulting coorbit spaces $\Coosp(\R^d \rtimes H, \Lsp_m^{p,q})$ related to classical function spaces like $\Bspq$ of $\Fspq$ (or their homogeneous counterparts)? Furthermore, how are coorbit spaces with respect to different dilation groups related to each other?

        This connects to the question posed above: If a given function/signal can be well approximated using one wavelet system, does the same also hold for a different system?
\end{enumerate}
Even for the special case of the shearlet dilation group, these questions are nontrivial and triggered several papers \cite{dakustte09,dastte12,dastte11,dahastte13}. Nevertheless, they also admit satisfactory answers in the present generality: 
As we will see, each of these questions is linked to the \textbf{dual action}
\[
    \varrho : H \times \R^d \to \R^d, (h,\xi) \mapsto h^{-T} \xi
\]
of the dilation group $H$ on the frequency space $\R^d$. To see the relevance of the dual action, note that the Fourier transform of $W_\varphi f(\cdot, h)$
is given by
\begin{equation}
    \left(\mathcal{F}[W_\varphi f (\cdot, h)]\right) (\xi) = |\det h|^{1/2} \cdot \widehat{f}(\xi) \cdot \overline{\widehat{\psi}(h^T \xi)}.
    \label{eq:WaveletTransformFourier}
\end{equation}
Thus, if $\widehat{\psi}$ has support in $U \subset \R^d$, then $W_\varphi f(\cdot, h)$ is bandlimited to $h^{-T}U = \varrho(h, U)$.

The remaining subsections deal with the three questions listed above.

\subsection{Question 1: Irreducibility and square-integrability of $\pi$}
\label{sub:DiscreteSeriesQuestion}
As shown in \cite{fu10,fu96,beta96}, the quasi-regular representation $\pi$ is irreducible and square-integrable if
and only if the following two properties are satisfied:
\begin{enumerate}
    \item There is some $\xi_0 \in \R^d \setminus \{0\}$ such that the orbit
          \[
              \calO := H^T \xi_0 = \left\{ h^T \xi_0 \with h \in H \right\}
          \]
          is \emph{open} and of full measure (i.e.\@ $\calO^c$ is a Lebesgue null set).

    \item The stabilizer
              $H_{\xi_0} := \left\{ h \in H \with h^T \xi_0 = \xi_0 \right\}$
          is compact.
\end{enumerate}
In this case, we call $\calO$ the (open) \textbf{dual orbit} of the dilation group $H$, whereas the null-set $\calO^c$ is called the \textbf{blind spot} of $H$.
Finally, a dilation group $H$ fulfilling the two properties   above is called \textbf{admissible}. In the following, we fix $\xi_0 \in \calO$.

To give the reader an idea of the richness of admissible dilation groups, we mention the following admissible dilation groups in dimension $d=2$:
\begin{enumerate}

    \item The \textbf{diagonal group}
        \[
            H_1 := \left\{ {\rm diag}(a,b) \with a,b \in \R \setminus \{0\} \right\}
        \]
        with dual orbit $\calO = \left( \R \setminus \{0\} \right)^2$.

    \item The \textbf{similitude group}
            $H_2 := (0,\infty) \cdot {\mathrm{SO}} (\R^2)$,
        with $\calO = \R^2 \setminus \{0\}$.

    \item The family of \textbf{shearlet type groups}
        \begin{equation}
            H_3^{(c)} := \left\{ \varepsilon \cdot \left( \begin{matrix}a & b \\ 0 & a^c \end{matrix} \right) \with a \in (0,\infty), b \in \R, \varepsilon \in \{1, -1\} \right\}.
            \label{eq:ShearletTypeGroup}
        \end{equation}
        Here, the \textbf{anisotropy parameter} $c \in \R$ can be chosen arbitrarily. Regardless of this choice, the dual orbit is always $\calO = (\R \setminus \{0\}) \times \R$.
\end{enumerate}

\subsection{Question 2: Existence of ``nice'' wavelets}
Here, we are given an admissible dilation group $H$ and we are interested in conditions which guarantee
$\calA_{v_0} \neq \left\{ 0 \right\}$ or $\calB_{v_0} \neq \{0\}$, where $v_0 : H \to (0,\infty)$
(interpreted by abuse of notation as a weight on $G = \R^d \rtimes H$) is a
(locally bounded) \textbf{control weight} for $\Ysp = \Lsp_m^{p,q}$. As we saw above, such a control weight always exists
under our general assumptions on the weight $m : H \to (0,\infty)$.
In the present setting, ``nice wavelets'' and ``better wavelets'' exist in abundance:
\begin{thm}
    \label{thm:BandlimitedWaveletsAreBetter}(cf.\@ \cite[Theorem 9]{fuvo15})

    Let $v_0 : H \to (0,\infty)$ be measurable and locally bounded. Then, every function $\varphi \in \mathcal{F}^{-1}(C_c^\infty (\calO)) \subset \ScRd$ satisfies $\varphi \in \calB_{v_0} \subset \calA_{v_0}$ and the map
 %
    \begin{equation}
        \varrho : C_c^\infty (\calO) \to \Lsp_{v_0}^{1}(G), g \mapsto W_\varphi (\mathcal{F}^{-1} g)
        \label{eq:BandLimitedWaveletsContinuousInclusionIntoReservoirPredual}
    \end{equation}
    is well-defined and continuous.
\end{thm}

This theorem, however, does \emph{not} yield existence of \emph{compactly supported} ``nice wavelets''.
In the case of ``traditional'' wavelets, it is well-known that a certain amount of
``vanishing moments'' (i.e.\@ $\partial^\alpha \widehat{\varphi}(0) = 0$ for all $|\alpha|\leq r$) makes a wavelet ``nice''.
This generalizes to the present setting:
\begin{defn}(\cite[Definition 1.4]{FuehrSimplifiedVanishingMomentCriteria})
    \label{defn:VanishingMoments}
    Let $r \in \N$. We say that $\varphi \in \Lsp^1(\R^d)$ has vanishing moments on $\calO^c$ of order $r$ if $\widehat{\varphi} \in C^r (\R^d)$ and if
    \[
        (\partial^\alpha \widehat{\varphi}) |_{\calO^c} \equiv 0 \qquad \text{ for } |\alpha| < r.
    \]
\end{defn}
As shown in \cite[Theorems 1.5(c), 2.12]{FuehrSimplifiedVanishingMomentCriteria}, given the weight $v_0 : H \to (0,\infty)$, one can explicitly
compute a number $\ell = \ell(H, v_0) \in \N$ such that every function $\psi \in \Lsp^1 (\R^d)$ with $\|\widehat{\psi}\|_{\ell, \ell}< \infty$
and with vanishing moments of order $\ell$ on $\calO^c$ satisfies $\psi \in \calB_{v_0} \subset \calA_{v_0}$.
Here, we employed the usual Schwartz-type norm
\[
    \|f\|_{\ell,\ell} := \max_{\alpha \in \N_0^d, |\alpha| \leq \ell} \,\, \sup_{x \in \R^d}  (1+|x|)^\ell \cdot |\partial^{\alpha} f(x)| \in [0,\infty].
\]
Finally, in the setting of ordinary wavelets, one can obtain a wavelet with enough vanishing moments by taking derivatives of a given smooth function
of sufficient decay.  In the present general setting, \cite[Lemmas 3.1 and 4.1]{fu12-1} yield a similar ``algorithm'' for obtaining functions
with suitably many vanishing moments on $\calO^c$. For the sake of brevity, we refrain from giving more details and instead refer the interested reader
to \cite{fu12-1, FuehrSimplifiedVanishingMomentCriteria}.

\subsection{Question 3: Relation of generalized wavelet coorbit spaces to other spaces}
\label{sub:CoorbitAsDecomposition}
Here, we are interested in the relation of the coorbit space $\Coosp(\Lsp_m^{p,q})$  to the classical
Besov- or Triebel-Lizorkin spaces, but also in the relation between $\Coosp(\Lsp_{m_1}^{p_1,q_1}, \R^d \rtimes H_1)$ and
$\Coosp(\Lsp_{m_2}^{p_2,q_2}, \R^d \rtimes H_2)$ for two different admissible dilation groups $H_1, H_2 \leq \GL (\R^d)$.

In view of the embedding theory for decomposition spaces (cf.\@ Subsection 
\ref{sub:DecompositionEmbedding}), this question would be solved (at least to a significant extent) if we knew that the coorbit space $\Coosp (\Lsp_m^{p,q})$
is (canonically isomorphic to) a decomposition space.
But the main result of \cite{fuvo15} is precisely such an isomorphism. In more detail, \cite{fuvo15} shows
\begin{equation}
    \Coosp(\Lsp_m^{p,q}(\R^d \rtimes H)) \cong \decompsp{\calQ_H}{p}{\lsp_u^q}
    \label{eq:CoorbitDecompositionIsomorphismAbstract}
\end{equation}
for a so-called \textbf{induced covering} $\calQ_H$ and a so-called \textbf{decomposition weight} $u$.
Below, we provide a concrete example  indicating  how this isomorphism
and the above embedding results  can be used to obtain
novel embedding results for shearlet coorbit spaces.

In the above isomorphism, the \textbf{induced covering} $\calQ_H$ of the dual orbit $\calO$ is given by $\calQ_H = (h_i^{-T}Q)_{i \in I}$
for certain $h_i \in H$ and a suitable $Q \subset \calO$. Its precise construction is  described in the following paragraph.
But once $\calQ = (h_i^{-T}Q)_{i \in I}$ is known, the \textbf{decomposition weight} $u = (u_i)_{i \in I}$ from above
(cf.\@ \cite[Lemma 35 and the ensuing remark]{fuvo15}) is given by
\begin{equation}
 u_i = |\det h_i|^{\frac{1}{2} - \frac{1}{q}} \cdot m(h_i) \qquad \forall i \in I.
    \label{eq:CoorbitDecompositionWeight}
\end{equation}

Let us reconsider the induced covering $\calQ_H$.  
The family $(h_i)_{i \in I}$ 
has to be \textbf{well spread} in $H$,
i.e.\@ 
there are compact unit neighborhoods $K_1, K_2 \subset H$ such that
\[
    H = \bigcup_{i \in I} h_i K_1 \qquad \text{ and } \qquad h_i K_2 \cap h_j K_2 = \emptyset \text{ for } i,j \in I \text{ with } i\neq j.
\]
Finally, the set $Q \subset \calO$ is an arbitrary open, bounded set such that the closure $\overline{Q} \subset \R^d$ is contained in $\calO$ and such that there
is a smaller open set $P \subset \R^d$ with $\overline{P} \subset Q$ and $\calO = \bigcup_{i \in I}h_i^{-T}P$. As shown in \cite[Theorem 20]{fuvo15}, such
sets $P, Q$ always exist if $(h_i)_{i \in I}$ is well-spread in $H$. The same theorem also shows that the resulting covering $\calQ_H = (h_i^{-T} Q)_{i \in I}$
is an (almost) structured admissible covering of $\calO$. In particular, there is a BAPU $(\varphi_i)_{i \in I}$ subordinate
to $\calQ_H$, cf.\@ Theorem \ref{thm:AlmostStructuredCoveringsAdmitBAPU}.

Finally, \cite[Lemmas 22 and 23]{fuvo15} show that the decomposition weight from equation
(\ref{eq:CoorbitDecompositionWeight}) is $\calQ_H$-moderate. This implies (cf.\@ Subsection \ref{sub:DecompositionSpaces}) that the decomposition space
$\decompsp{\calQ_H}{p}{\lsp_u^q}$ on the right-hand side of equation (\ref{eq:CoorbitDecompositionIsomorphismAbstract}) is well-defined.

\begin{figure}[h]
    \centering
    \includegraphics[scale=0.27]{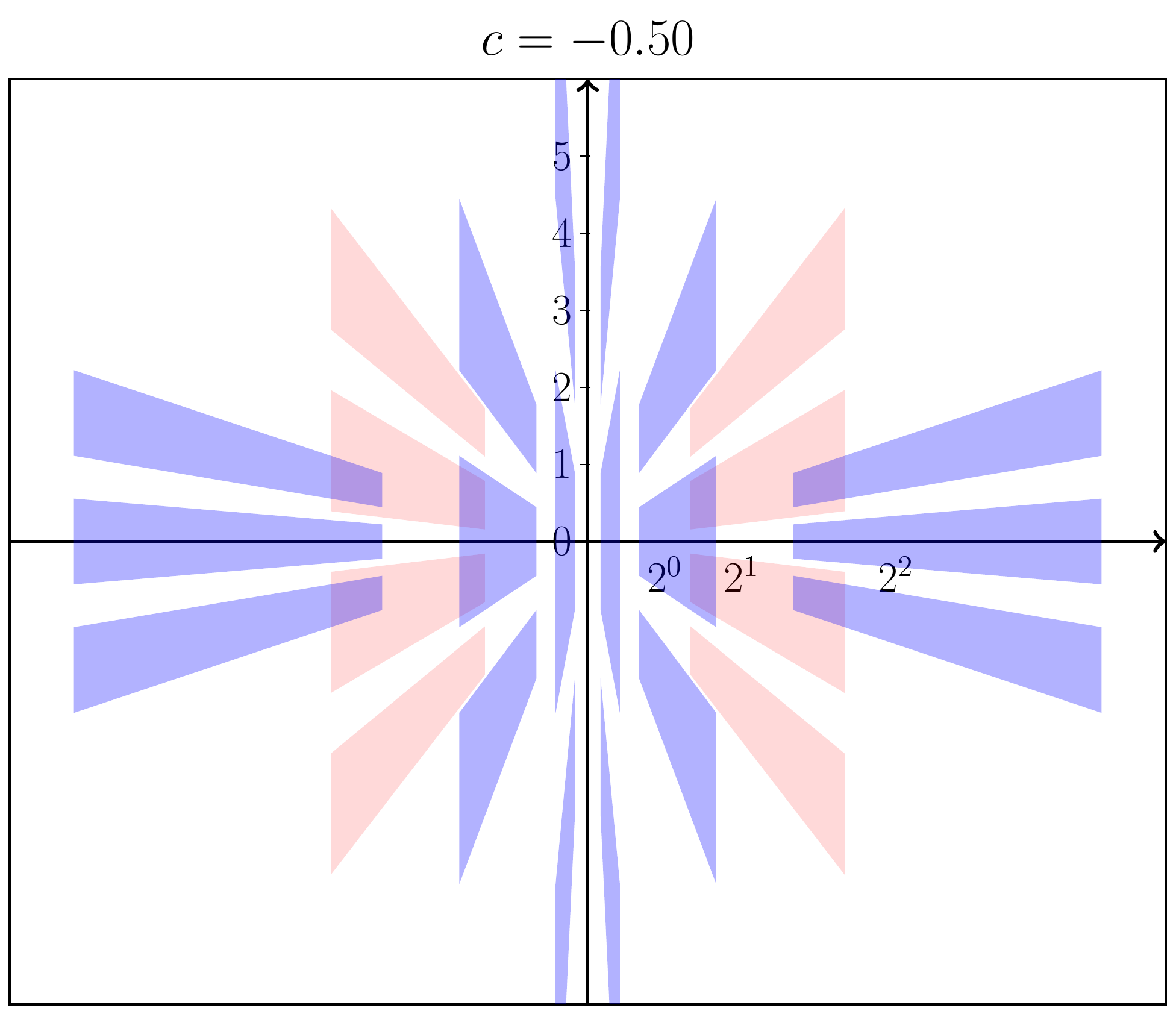}
    \includegraphics[scale=0.27]{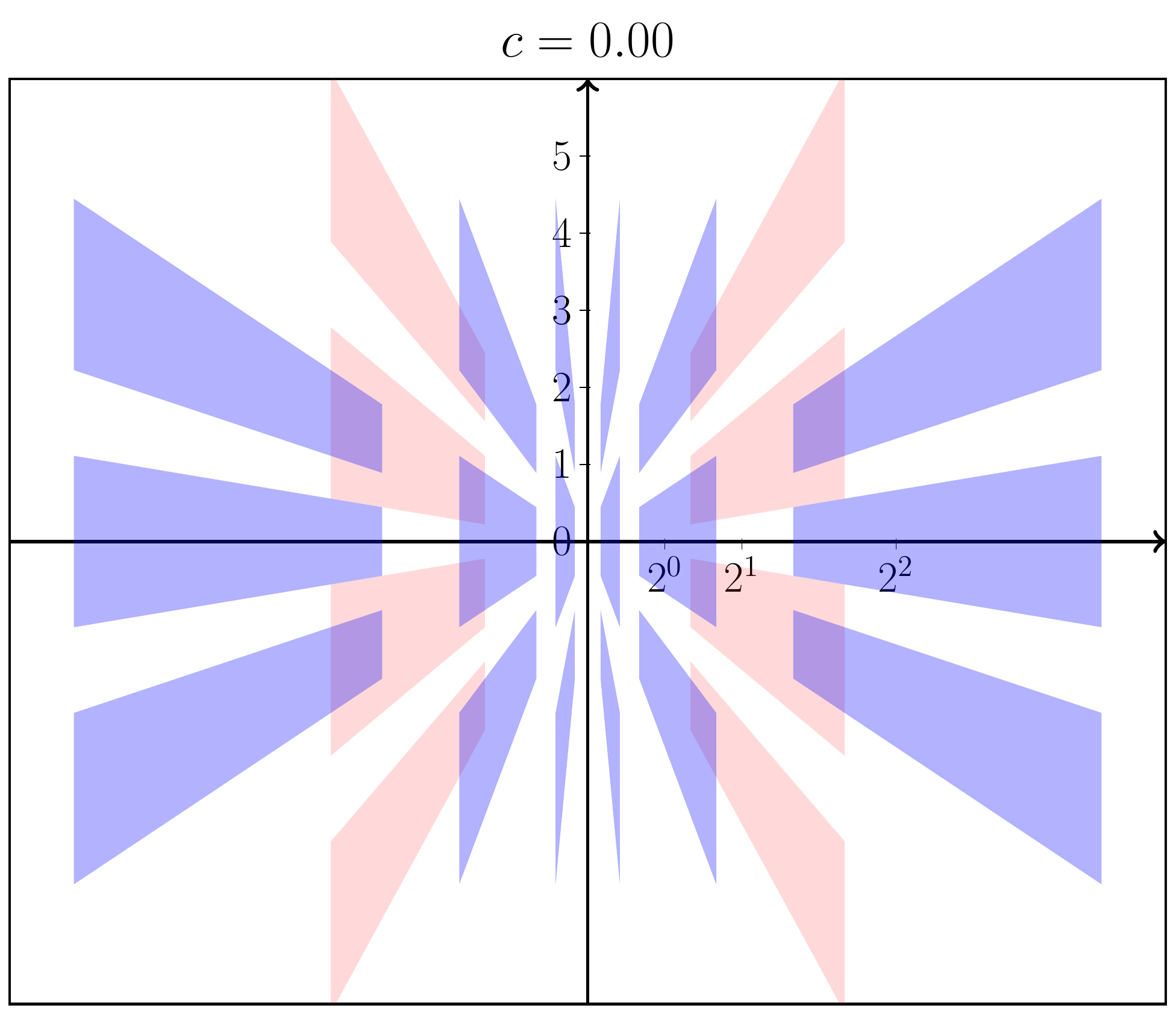}\\
    \includegraphics[scale=0.27]{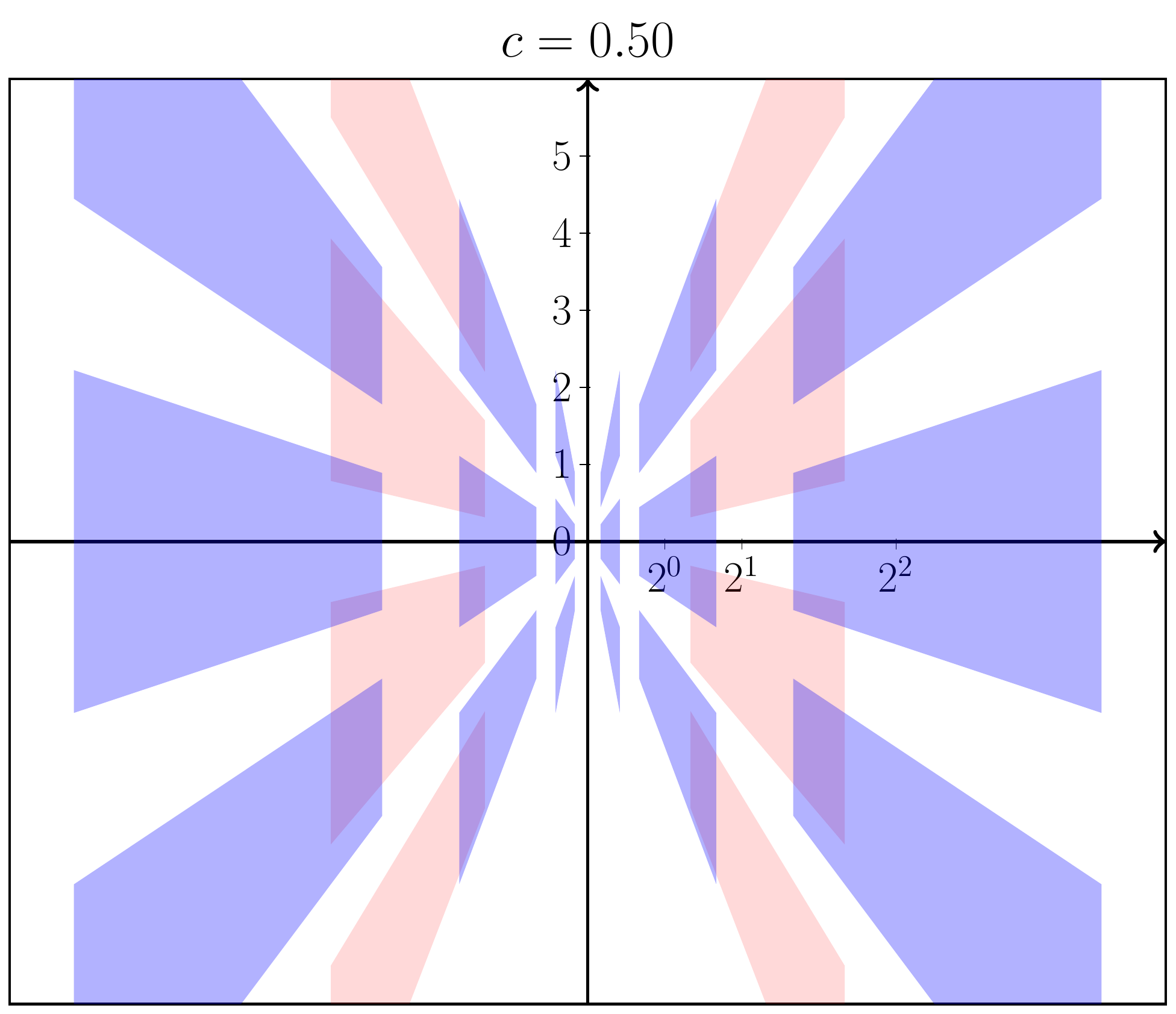}
    \caption{The figure shows (a part of) the induced covering for the group $H_3^{(c)}$ for $c= - \frac{1}{2}$, $c=0$ and $c = \frac{1}{2}$. These choices show the qualitatively different behaviour of the covering for different values of $c$. For $c_1 < c_2$, the covering $\calS^{(c_1)}$ is ``larger/coarser'' near the $y$-axis, whereas the covering $\calS^{(c_2)}$ is ``larger/coarser'' away from the $y$-axis.
    Note: The images are taken from \cite[Figure 1]{VoigtlaenderPhDThesis}.}
    \label{fig:ShearletCovering}
\end{figure}

As an illustration of the concept of an \emph{induced covering}, we consider the shearlet type group $H_3^{(c)}$ from equation (\ref{eq:ShearletTypeGroup}).
A picture of (a part of) the associated induced covering for different values of the \emph{anisotropy parameter}
$c \in \R$ is shown in Fig.\@ \ref{fig:ShearletCovering}.

For the explicit description of the isomorphism from equation (\ref{eq:CoorbitDecompositionIsomorphismAbstract}), we
note that it is a (relatively) easy consequence of equation (\ref{eq:BandLimitedWaveletsContinuousInclusionIntoReservoirPredual})
that
\[
    \iota : Z(\calO) = \mathcal{F}(C_c^\infty (\calO)) \to \Hilb^{1}_{v_0}, f \mapsto \overline{f},
\]
with $\Hilb_{v_0}^{1}$ as in equation (\ref{eq:CoorbitTheoryTestFunctions}) is well-defined, antilinear and continuous.
By duality, this implies that
\begin{equation}
    \iota^T : \left( \Hilb_{v_0}^1 \right)^\angle \to Z'(\calO), \theta \mapsto \theta \circ \iota
    \label{eq:CoorbitIsomorphism}
\end{equation}
is well-defined, continuous and \emph{linear}. Recall from Section \ref{sec:CoorbitTheory} that $\Coosp(\Lsp_m^{p,q})$ is a subspace of the \emph{reservoir}
$\calR = \left( \Hilb_{v_0}^1 \right)^\angle$ for a certain control weight $v_0$.
Thus, $\iota^T f \in Z'(\calO)$ is well-defined for every $f \in \Coosp(\Lsp_m^{p,q})$. The claim of equation (\ref{eq:CoorbitDecompositionIsomorphismAbstract})
is precisely that $\iota^T : \Coosp(\Lsp_m^{p,q}) \to \decompsp{\calQ}{p}{\lsp_u^q}$ is an isomorphism of Banach spaces. A proof of this
fact can be found in \cite[Theorem 43]{fuvo15}.

This representation of $\Coosp(\Lsp_m^{p,q}(G))$ as a decomposition space has several important consequences, both mathematically and conceptually:
\begin{itemize}
    \item As we saw above, the coorbit space $\Coosp(\Lsp_m^{p,q}(G))$ with its original definition is heavily tied to the group
        $G = \R^d \times H$ and thus to the dilation group $H$. In particular, $\Coosp(\Lsp_m^{p,q}(G))$ is a subspace of the reservoir
        $\calR = (\calH_{v_0}^1 (G))^\angle$.

        This makes it difficult to consider an element $f \in \Coosp(\Lsp_{m_1}^{p_1,q_1}(\R^d \rtimes H_1))$ as an element of another coorbit space
        $\Coosp(\Lsp_{m_2}^{p_2,q_2}(\R^d \rtimes H_2))$, or as an element of more classical function spaces like $\Bspq$.

        In contrast, as we saw in Section \ref{sec:DecompositionSpaces}, it is (at least in principle) possible
        to compare decomposition spaces $\decompsp{\calQ}{p_1}{\lsp_u^{q_1}}$
        and $\decompsp{\calP}{p_2}{\lsp_v^{q_2}}$ which are defined using two different coverings $\calQ, \calP$ of the sets $\calO, \calO' \subset \R^d$.

        Thus, using the decomposition space view, it becomes possible to compare wavelet coorbit spaces defined by  \emph{different} dilation groups.

    \item Even ignoring the issue of the different reservoirs (e.g.\@ by restricting to Schwartz functions), it is not at all obvious how the decay
        or integrability condition $W_\varphi f \in \Lsp_{m_1}^{p_1,q_1}(\R^d \rtimes H_1)$ relates to another decay condition
        $W_\varphi f \in \Lsp_{m_2}^{p_2,q_2}(\R^d \rtimes H_2)$, even if the \emph{same} analyzing window is used in both cases. One of the reasons is that it is difficult to compare the two actions of the dilation groups on $\varphi$, as well as the two distinct Haar measures.

        In comparison, the decomposition space point of view translates these two elusive properties into (more or less) transparent quantities, namely
        \begin{enumerate}
            \item The \textbf{induced covering} $\calQ_H = (h_i^{-T} Q)_{i \in I}$ for some well-spread family $(h_i)_{i \in I}$ in $H$
                  and a suitable set $Q \subset \calO$,
            \item The \textbf{decomposition weight}
            $u_i = |\det h_i|^{\frac{1}{2} - \frac{1}{q}} \cdot m(h_i)$.
        \end{enumerate}

        Using the methods from Subsection \ref{sub:DecompositionEmbedding}, it is then (comparatively) easy to establish embeddings
            $\decompsp{\calQ_{H_1}}{p_1}{\lsp_{u_1}^{q_1}} \hookrightarrow \decompsp{\calQ_{H_2}}{p_2}{\lsp_{u_2}^{q_2}}$
        between the associated decomposition spaces
        and thus of the two coorbit spaces $\Coosp(\Lsp_{m_1}^{p_1, q_1}(\R^d \rtimes H_1))$ and $\Coosp(\Lsp_{m_2}^{p_2, q_2}(\R^d \rtimes H_2))$.

    \item Similarly, one can use the methods from Subsection \ref{sub:DecompositionEmbedding} to establish embeddings between generalized wavelet coorbit spaces and classical smoothness spaces like Sobolev- and Besov spaces.

    \item Conceptually, all these considerations show that \emph{the approximation theoretic properties of the wavelet system generated by a dilation group $H$ are completely determined by the way in which (the dual action of) $H$ covers/partitions the frequencies.}
\end{itemize}

\noindent
Theorem \ref{thm:ShearletIntoBesov} below is an example of results that can be obtained by combining
the embedding results from Subsection 
\ref{sub:DecompositionEmbedding}
with the isomorphism
\[
    \Coosp(\Lsp_m^{p,q}(\R^d \rtimes H)) \cong \decompsp{\calQ_H}{p}{\lsp_u^q}.
\]
Precisely, we consider embeddings between shearlet coorbit spaces and {\it inhomogeneous} Besov spaces. For the sake of brevity, we only consider embeddings
of the shearlet coorbit space into inhomogeneous Besov spaces.
Results for the reverse direction are also available (cf.\@ \cite[Theorem 6.3.14]{VoigtlaenderPhDThesis}), but are omitted here.

We only consider the case $c \in (0,1]$. This ensures that the induced covering $\calS^{(c)} = \calQ_{H_3^{(c)}}$ is
almost subordinate to the inhomogeneous dyadic covering. See \cite[Lemma 6.3.10]{VoigtlaenderPhDThesis} for a formal
proof and Figure \ref{fig:ShearletCovering} for a graphical illustration.
Note, however, that $\calS^{(c)}$ is \emph{not} relatively moderate with respect to the inhomogeneous dyadic covering.
This limits sharpness of our results to a certain range of $p_2$.
\begin{thm}(\cite[Theorem 6.3.14]{VoigtlaenderPhDThesis})
    \label{thm:ShearletIntoBesov}
Let $c \in (0,1]$, $p_1, p_2, q_1, q_2 \in [1,\infty]$ and $\alpha, \beta, \gamma \in \R$. Set $p_2^{\triangledown} := \min \{p_2, p_2 '\}$, define the weight
    \[
        u^{(\alpha, \beta)} : H_3^{(c)} \to (0,\infty), h \mapsto \|h^{-1}\|^{\alpha} \cdot |\det h|^{\beta}
    \]
    and set
    \begin{align*}
        \alpha^{(1)} &:= \frac{1+c}{c} \left( \frac{1}{p_1} - \frac{1}{p_2} - \frac{1}{q_1} + \frac{1}{2} + \beta \right), \\
        \gamma^{(1)} &:= - (1+c) \cdot \left( \frac{1}{p_1} - \frac{1}{p_2} - \frac{1}{q_1} + \frac{1}{2} + \beta \right) + (c-1) \cdot \left( \frac{1}{p_2^{\triangledown}} - \frac{1}{q_1} \right)_{+}.
    \end{align*}

    If $p_1 \leq p_2$ as well as
    \[
        \begin{cases}
            \gamma \leq \alpha + \gamma^{(1)}, &\text{if } q_1 \leq q_2, \\
            \gamma <    \alpha + \gamma^{(1)}, &\text{if } q_1 >    q_2
        \end{cases}
        \quad \text{ and } \quad
        \begin{cases}
            \max \left\{ \frac{1}{p_2^\triangledown} - \frac{1}{q_1}, \, \alpha \right\} < \alpha^{(1)}, & \text{if } q_1 > p_2^{\triangledown}, \\
            \max \left\{ 0,\alpha \right\} \leq \alpha^{(1)},                                            & \text{if } q_1 \leq p_2^{\triangledown}
        \end{cases}
    \]
    hold, then
    \begin{equation}
        \Coosp(\Lsp_{u^{(\alpha, \beta)}}^{p_1, q_1}(\R^2 \rtimes H_3^{(c)})) \hookrightarrow \Bsp^{\gamma}_{p_2, q_2} (\R^2).
        \label{eq:ShearletCoorbitIntoBesov}
    \end{equation}
    A necessary condition for existence of this embedding is obtained by replacing $p_2^{\triangledown}$ by $p_2$ everywhere (also in the definition of $\gamma^{(1)}$).\\[0.2cm]
\end{thm}
\begin{rem*}
    \begin{enumerate}
        \item Existence of the embedding (\ref{eq:ShearletCoorbitIntoBesov}) has to be interpreted suitably.
            Precisely, (\ref{eq:ShearletCoorbitIntoBesov}) means that there is a bounded linear map
            \[
                \iota : \Coosp(\Lsp_{u^{(\alpha, \beta)}}^{p_1, q_1}(\R^2 \rtimes H_3^{(c)})) \to \Bsp_{p_2, q_2}^{\gamma}(\R^2)
            \]
            which satisfies $\iota f = f$ for all $f \in \Lsp^2 (\R^2) \cap \Coosp(\Lsp_{u^{(\alpha, \beta)}}^{p_1, q_1}(\R^2 \rtimes H_3^{(c)}))$.

        \item The preceding theorem is superficially similar to \cite[Theorem 4.7]{dastte11}. But the two results are very different,
            since Dahlke et al.\@ consider embeddings of the \emph{strict subspace} $\mathcal{SCC}_{p,r} \lneq \Coosp(\Lsp_{u^{(0,-2r/3)}}^{p,p}(\R^2 \rtimes H_3^{(1/2)}))$ into a \emph{sum} of \emph{homogeneous} Besov spaces $\dot{\Bsp}_{p,p}^{\sigma_1}(\R^2) + \dot{\Bsp}_{p,p}^{\sigma_2}(\R^2)$ for certain $\sigma_1, \sigma_2$.

            In contrast, the preceding theorem investigates embeddings of the \emph{whole} shearlet coorbit space into a \emph{single}, \emph{in}homogeneous Besov space.

        \item The preceding theorem achieves a \emph{complete characterization} of the embedding (\ref{eq:ShearletCoorbitIntoBesov}) for $p_2 \in [1,2]$, since we have $p_2^{\triangledown} = p_2$ in this range.
    \end{enumerate}
\end{rem*}

As a conclusion, we remark that the embedding results and the isomorphism between generalized wavelet coorbit spaces and decomposition spaces can also be used
to derive embeddings between the coorbit spaces $\Coosp(\Lsp_m^{p,q}(\R^2 \rtimes H_3^{(c)}))$ for different values of $c$, see \cite[Theorem 6.3.9]{VoigtlaenderPhDThesis}.

They can also be used to derive (non)boundedness of certain operators---e.g.\@ dilation operators---acting on coorbit spaces.
For example, in \cite[Theorem 6.5.9]{VoigtlaenderPhDThesis}, the set of matrices which act boundedly by dilation simultaneously on \emph{all}
coorbit spaces of the shearlet type group $H_3^{(c)}$ (for a fixed $c \in (0,1)$) is characterized completely.

%

\section*{Acknowledgments}
Both authors want to thank HIM---the Hausdorff Institute of
Mathematics---where we both spent some time during the preparation
of this manuscript.
FV was funded by the Excellence Initiative of the German federal and state governments,
and by the German Research Foundation (DFG),
under the contract FU 402/5-1.


\bibliographystyle{abbrv}


\end{document}